\documentclass[12pt,leqno]{amsart}
\usepackage{amsmath,amsthm,amsfonts,amsbsy,multicol,fancybox}
\usepackage[utf8]{inputenc}
\usepackage[T1]{fontenc}

\usepackage{caption}\usepackage{subcaption}
\usepackage{graphicx}\usepackage{amssymb,latexsym}
\usepackage{epstopdf}
\usepackage{color}
\DeclareGraphicsExtensions{.pdf, .jpg, .pnp, .bmp, .fig}
\usepackage{tikz}
\usepackage[all]{xy}
\input xy
\xyoption{all}

\usepackage[hidelinks]{hyperref}

\usepackage{thmtools}

\newcounter{argument}
\newenvironment{argument}[1][\medskip]{%
\refstepcounter{argument}
\par\medskip
\noindent\phantomsection
\textbf{#1~\thesection.\arabic{argument}\,\,}\rmfamily\em}{\hspace{\fill}$\Box$\par\smallskip\noindent}

\newcommand{\bass}{\begin{argument}[Assumption]}\newcommand{\eass}{\end{argument}}
\newcommand{\bth}{\begin{argument}[Theorem]} \newcommand{\ethe}{\end{argument}}
\newcommand{\bre}{\begin{argument}[Remark]}      \newcommand{\ere}{\end{argument}}
\newcommand{\ble}{\begin{argument}[Lemma]}       \newcommand{\ele}{\end{argument}}
\newcommand{\bde}{\begin{argument}[Definition]}   \newcommand{\ede}{\end{argument}}
\newcommand{\bco}{\begin{argument}[Corollary]}     \newcommand{\eco}{\end{argument}}
\newcommand{\bpr}{\begin{argument}[Proposition]}  \newcommand{\epr}{\end{argument}}
\newcommand{\bexam}{\begin{argument}[Example]}\newcommand{\eexam}{\end{argument}}

\newcommand{\bpf}{\begin{proof}}\newcommand{\epf}{\end{proof}}
\newcommand{\barr}{\begin{array}}\newcommand{\earr}{\end{array}}
\newcommand{\beao}{\begin{eqnarray*}}\newcommand{\eeao}{\end{eqnarray*}\noindent}
\newcommand{\beam}{\begin{eqnarray}}\newcommand{\eeam}{\end{eqnarray}\noindent}
\newcommand{\beqq}{\begin{equation}}\newcommand{\eeqq}{\end{equation}\noindent}

\newcommand{\ov}{\overline} \newcommand{\un}{\underbrace}

\newcommand{\D}{\Delta}
  \newcommand{\ep}{\epsilon}

\newcommand{\lam}{\lambda} 

\newcommand{\w}{\omega}

 \newcommand{\bbC}{{\mathbb C}}

\newcommand{\bfE}{{\mathbb E}}

\newcommand{\bbi}{{\mathbb I}}

 \newcommand{\bbN}{{\mathbb N}}

\newcommand{\bbr}{{\mathcal R}} \newcommand{\bbR}{{\mathbb R}}
\newcommand{\bbs}{{\mathcal S}}

\graphicspath{{./../Images/}}

\begin{document}

\title[A note on Asymptotic mean-square stability for SDEs]{A note on Asymptotic mean-square stability of stochastic linear two-step methods for SDEs}

\author[I. S. Stamatiou]{I. S. Stamatiou}
\email{joniou@gmail.com}

\begin{abstract}
In this note we study the asymptotic mean-square stability for two-step schemes applied to a  scalar stochastic differential equation (sde) and applied to systems of sdes. We derive necessary and sufficient conditions for the asymptotic MS-stability of the methods in terms of the parameters of the schemes. The stochastic Backward Differentiation Formula (BDF2) scheme is asymptotically mean-square stable for any step-size whereas the two-step Adams-Bashforth (AB2) and Adams-Moulton (AM2) methods are unconditionally stable. The improved versions of the schemes do not perform better w.r.t their stability behavior in the scalar case, as expected, but  the situation is different in more dimensions. Numerical experiments confirm theoretical results.
\end{abstract}
\date\today

\maketitle

\tableofcontents

\section{Introduction.}\label{NSF:sec:intro}
\setcounter{equation}{0}

Consider the general type $d$-dimensional It\^o stochastic differential equation  (sde) 
\beqq\label{STB-eq:sde.general}
dX(t) = F(t, X(t)) dt + G(t, X(t)) dW(t), \quad X(t_0)=X_0,
\eeqq
driven by the $m$-dimensional Wiener process, where the coefficients $F: [0,T]\times\bbR^d\mapsto\bbR^d, G: [0,T]\times\bbR^d\mapsto\bbR^{d\times m}$ are such that there exists a unique path-wise strong solution of (\ref{STB-eq:sde.general}), cf. \cite[Ch. 2.3]{mao:2007}. We will also study complex-valued functions $F, G, X.$ 
The two-step Maryuama method with an equidistant step-size $h$ for the approximations $X_i\approx X(t_i)$ of the solution of (\ref{STB-eq:sde.general}) read
\beqq\label{STB-eq:twostepapp_gen}
\sum_{j=0}^{2}\alpha_jX_{i-j} = h\sum_{j=0}^{2} \beta_jF_{i-j} + \sum_{j=1}^{2}\sum_{r=1}^{m} \gamma_jG_{r,i-j}\sqrt{h}\xi_{r,i-j}, \quad i=2, 3, \ldots,
\eeqq
and the improved two-step Maryuama method is given by 
\beqq\label{STB-eq:imptwostepapp_gen}
\sum_{j=0}^{2}\alpha_jX_{i-j} = h\sum_{j=0}^{2} \beta_jF_{i-j} + \sum_{j=1}^{2}\sum_{r=1}^{m} \left(\gamma_jG_{r,i-j}\sqrt{h}\xi_{r,i-j} + (\gamma_j+\eta_j)(F^{\prime}G_{r})_{i-j}h^{3/2}\xi_{r,i-j}\right), 
\eeqq
for $i=2,3,\ldots,$ in the case of systems with commutative noise; here $\{\xi_{r,i}\}_{i\in\bbN_0}, r=1, \ldots, m,$ are sequences of i.i.d. standard normal r.v.s, $(\alpha_j, \beta_j), j=0,\ldots,2$ and $(\gamma_j, \eta_j), j=1,2$ are appropriate parameters and $f_{i-j}$ denotes  $f(t_{i-j},X_{i-j})$ for appropriate functions $f$ as above. For convergence properties of (\ref{STB-eq:twostepapp_gen}) and (\ref{STB-eq:imptwostepapp_gen}) see \cite{buckwar_winkler:06}, \cite{buckwar_winkler:07}. Here, we are interested in  mean-square asymptotic properties of the above numerical approximations. We perform a linear stability analysis using linear time-invariant test equations; in \cite{buckwar_etal:06} sufficient conditions are given for asymptotic mean-square stability  of (\ref{STB-eq:twostepapp_gen}) applying appropriate Lyapunov-type functionals. We provide in our main result, Theorem \ref{STB-thr:2step}, necessary and sufficient conditions  for the asymptotic mean-square stability  of (\ref{STB-eq:twostepapp_gen}) and (\ref{STB-eq:imptwostepapp_gen}) following a different approach. 

In Section \ref{NSF:sec:lineartesteqn} we use the scalar linear test-equation to study the stability properties of the two-step Maruyama methods. The stability matrix $\bbs$ of the two-step methods is analyzed.  Section \ref{NSF:sec:conditions} provides our main result regarding  stability conditions for two-step Maruyama methods and applications of it. The linear mean-square stability of the methods is studied in Section \ref{NSF:sec:linearMSstability} and experiments are made in Section \ref{NSF:sec:exp}. Section \ref{NSF:sec:system} is devoted to systems of linear test equations with multi-dimensional noise.

\section{Linear Test equation.}\label{NSF:sec:lineartesteqn}
\setcounter{equation}{0}

Consider the scalar linear test-equation with multiplicative noise
\beqq\label{STB-eq:sde.scalar}
dX(t) = \lam X(t) dt + \mu X(t) dW(t), \quad X(t_0)=X_0,
\eeqq
where the coefficients $\lam,\mu\in\bbC$ and assume w.l.o.g. that $X_0$ is non-random.
The two-step Maryuama method with an equidistant step-size $h$ for the approximations $X_i\approx X(t_i)$ of the solution of (\ref{STB-eq:sde.scalar}) read, (apply (\ref{STB-eq:twostepapp_gen}) with $F_i=\lam X_i, G_i=\mu X_i, m=1$)
\beqq\label{STB-eq:twostepapp}
\sum_{j=0}^{2}\alpha_jX_{i-j} = h\sum_{j=0}^{2} \beta_j\lam X_{i-j} + \sum_{j=1}^{2} \gamma_j\mu X_{i-j}\sqrt{h}\xi_{i-j}, \quad i=2, 3, \ldots,
\eeqq
and the improved two-step Maryuama method is given by 
\beqq\label{STB-eq:imptwostepapp}
\sum_{j=0}^{2}\alpha_jX_{i-j} = h\sum_{j=0}^{2} \beta_j\lam X_{i-j} + \sum_{j=1}^{2} \left(\gamma_j\mu X_{i-j}\sqrt{h}\xi_{i-j} + (\gamma_j+\eta_j)\lam\mu X_{i-j}h^{3/2}\xi_{i-j}\right), i=2,3,\ldots,
\eeqq
where $\{\xi_i\}_{i\in\bbN_0}$ is a sequence of i.i.d. standard normal r.v.s and $(\alpha_j, \beta_j), j=0,\ldots,2$ and $(\gamma_j, \eta_j), j=1,2$ are appropriate parameters.

The recurrences (\ref{STB-eq:twostepapp}) and (\ref{STB-eq:imptwostepapp}) can be rewritten in the form
\beqq\label{STB-eq:scalar.diff}
X_{i} = aX_{i-1} + cX_{i-2} + bX_{i-1}\xi_{i-1} + dX_{i-2}\xi_{i-2}, \, i=2,3,\ldots,
\eeqq
where for (\ref{STB-eq:twostepapp}) the complex coefficients $a,b,c$ and $d$ read 
\beqq\label{STB-eq:twostepappcoef}
a=\frac{-\alpha_1 + h\beta_1\lam}{\alpha_0-h\beta_0\lam}, \quad b=\frac{\sqrt{h}\gamma_1\mu}{\alpha_0-h\beta_0\lam}, \quad c=\frac{-\alpha_2 + h\beta_2\lam}{\alpha_0-h\beta_0\lam}, \quad d=\frac{\sqrt{h}\gamma_2\mu}{\alpha_0-h\beta_0\lam}
\eeqq
and for (\ref{STB-eq:imptwostepapp}) the complex coefficients $a,c$ are the same and $b, d$ read
\beqq\label{STB-eq:twostepapp1coef}
b^*= b + \frac{\lam\mu(\gamma_1+\eta_1)h^{3/2}}{\alpha_0-h\beta_0\lam}, \quad d^* = d + \frac{\lam\mu(\gamma_2+\eta_2)h^{3/2} }{\alpha_0-h\beta_0\lam}.
\eeqq
In Table \ref{STB-tab:coefs} we list the coefficients $\alpha_i, \beta_i, \gamma_i, \eta_i$ for different two-step Maruyama methods.
\begin{table}\caption{Coefficients of two-step Maruyama schemes and improved two-step schemes as in (\ref{STB-eq:twostepapp}) and (\ref{STB-eq:imptwostepapp}) with $\alpha_0=\gamma_1=1.$}\label{STB-tab:coefs}
 \begin{tabular}{r|| r r r r r r r r} 
  \hline
Method & $\alpha_1$ & $\alpha_2$ & $\beta_0$ & $\beta_1$ & $\beta_2$ & $\gamma_2$ & $\eta_1$ & $\eta_2$\\ [0.5ex] 
 \hline\hline
 AB2 & $-1$ & $0$ & $0$ & $3/2$ & -$1/2$ &  $0$ & $-$ & $-$\\ [0.5ex] \hline
 AB2I & $-1$ & $0$ & $0$ & $3/2$ & -$1/2$ &  $0$ & $0$ & $-1/2$\\ [0.5ex]\hline
 AM2 & $-1$ & $0$ & $5/12$ & $8/12$ & -$1/12$ &  $0$ & $-$ & $-$\\ [0.5ex] \hline
 AM2I & $-1$ & $0$ & $5/12$ & $8/12$ & -$1/12$ &  $0$ & $-5/12$ & $-1/12$\\ [0.5ex]\hline
 BDF2 & $-4/3$ & $1/3$ & $2/3$ & $0$ & $0$ & $-1/3$ & $-$ & $-$\\ [0.5ex]\hline
BDF2I & $-4/3$ & $1/3$ & $2/3$ & $0$ & $0$ & $-1/3$ & $-2/3$ & $1/3$
\end{tabular}
\end{table}

 The stability or transition  matrix $\bbs$ of the two-step method (\ref{STB-eq:scalar.diff}) reads
\beqq\label{STB-eq:Smat.scalar}
\bbs= \begin{bmatrix}
    |a|^2 + |b|^2 &  \ov{a}c & a\ov{c} & |c|^2 + |d|^2 + ab\ov{d} +\ov{a}\ov{b}d \\
    \ov{a} & 0 & \ov{c} & b\ov{d} \\
    a & c & 0 & \ov{b}d \\
    1 & 0 & 0  & 0
\end{bmatrix},
\eeqq
where $\ov{z}$ stands for the conjugate of $z\in\bbC.$ The zero solution of the difference equations (\ref{STB-eq:scalar.diff}) is asymptotically mean-square stable iff the spectral radius of the mean-square stability matrix $\bbs$ satisfies
\beqq\label{STB-eq:stabcond}
\rho(\bbs)<1.
\eeqq
Recall that $\rho(\bbs):=\max|l_j|$ where $l_j$ are the eigenvalues of $\bbs.$ Computing the eigenvalues of $\bbs$ amounts to finding the roots of its characteristic polynomial and verifying condition (\ref{STB-eq:stabcond}). Here the characteristic polynomial  is a fourth-order polynomial given by
\beqq\label{STB-eq:charpolSmat.scalar}
P(z) = z^4 + p_1z^3+p_2z^2 + p_3z + p_4
\eeqq       
 where the real coefficients $p_j,  j=1, \ldots, 4,$ read
\beqq\label{STB-eq:charpolcoef1Smat.scalar}
p_1=-|a|^2 - |b|^2, \quad p_2=-2|c|^2 - |d|^2 - 2\Re(ab\ov{d}) - 2\Re(a^2\ov{c}),
\eeqq
\beqq\label{STB-eq:charpolcoef2Smat.scalar}
p_3=-2\Re(\ov{a}bc\ov{d}) - |a|^2|c|^2 + |b|^2|c|^2, \quad p_4=|c|^4 + |c|^2|d|^2.
\eeqq
We can check $\rho(\bbs)<1$ avoiding the computation of the $\rho(\bbs)$ by verifying conditions on the parameters $p_j, j=1, \ldots, 4,$ implied by the Schur-Cohn criterion.   The strategy is the following (cf. \cite{jury:1988}): 
Define the transpose $P^{\#}$ of $P$ as 
$$
P^{\#}(z) = z^4\ov{P}\left(\frac{1}{\ov{z}}\right) = p_4z^4 + p_3z^3 + p_2z^2 + p_1z + 1;
$$
define the $4\times4$ Schur-Cohn matrix associated to $P$ by       
\beao
&&\D_4(P,P^{\#}) = \begin{bmatrix}
    1 &  0 & 0 & 0 \\
    p_1 & 1 & 0 & 0 \\
    p_2 & p_1&  1  & 0\\
    p_3 & p_2 & p_1  & 1
\end{bmatrix}\begin{bmatrix}
    1 &  0 & 0 & 0 \\
    p_1 & 1 & 0 & 0 \\
    p_2 & p_1&  1  & 0\\
    p_3 & p_2 & p_1  & 1
\end{bmatrix}^T - \begin{bmatrix}
    p_4 &  0 & 0 & 0 \\
    p_3 & p_4 & 0 & 0 \\
    p_2 & p_3 &  p_4  & 0\\
    p_1 & p_2 & p_3  & p_4
\end{bmatrix}\begin{bmatrix}
    p_4 &  0 & 0 & 0 \\
    p_3 & p_4 & 0 & 0 \\
    p_2 & p_3 &  p_4  & 0\\
    p_1 & p_2 & p_3  & p_4
\end{bmatrix}^T   \\
&=&\begin{bmatrix}
    1 &  0 & 0 & 0 \\
    p_1 & 1 & 0 & 0 \\
    p_2 & p_1&  1  & 0\\
    p_3 & p_2 & p_1  & 1
\end{bmatrix}\begin{bmatrix}
    1 &  p_1 & p_2 &p_3 \\
    0 & 1 & p_1 & p_2 \\
    0 & 0 &  1  & p_1\\
    0 & 0 & 0  & 1
\end{bmatrix} - \begin{bmatrix}
    p_4 &  0 & 0 & 0 \\
    p_3 & p_4 & 0 & 0 \\
    p_2 & p_3 &  p_4  & 0\\
    p_1 & p_2 & p_3  & p_4
\end{bmatrix}\begin{bmatrix}
    p_4 &  p_3 & p_2 & p_1 \\
    0 & p_4 & p_3 & p_2 \\
    0 & 0 &  p_4  & p_3\\
    0 & 0 & 0  & p_4
\end{bmatrix}\\
&=&\begin{bmatrix}
    1 - |p_4|^2 &  p_1 -  p_4p_3 & p_2 - p_4p_2&p_3 - p_4p_1 \\
    p_1 -  p_3p_4 & 1 + |p_1|^2 - |p_3|^2 - |p_4|^2 & p_1p_2 + p_1 - p_3p_2 - p_4p_3 & p_2  - p_4p_2  \\
    p_2  - p_2p_4& p_2p_1 + p_1 - p_2p_3 - p_3p_4&  1 + |p_1|^2 -|p_3|^2 - |p_4|^2& p_1  - p_4p_3\\
    p_3 - p_1p_4& p_2 - p_2p_4& p_1 - p_3p_4 & 1 - |p_4|^2 \end{bmatrix} 
\eeao
where $Q^T$ is the  transpose of a matrix $Q,$ i.e. the matrix with entries $Q^T_{ij}=Q_{ji}$;      
the polynomial P has all the roots inside the unit disk iff  $ \D_4(P,P^{\#})$ is positive definite, which corresponds to 
\beqq\label{STB-eq:SCcriterion.scalar}
\textup{det} \D_k(P,P^{\#})>0, \quad k=1, \ldots, 4.
\eeqq
In order to decide about  (\ref{STB-eq:SCcriterion.scalar}) we can use the connection it has with the Schur coefficients $(\nu_k)_{k=0,\ldots,3}$ of the pair $(P,P^{\#});$ in general for the pair $(P,Q)$ we construct the sequence $(P_k, Q_k)_{k=0,1,\ldots,}$ as
$P_0=P,Q_0=Q$ and 
$$
P_{k}(z) = \frac{1}{z}\left(P_{k-1}(z) - \nu_{k-1}Q_{k-1}(z)\right), \quad k\geq1,
$$
$$
Q_{k}(z) = Q_{k-1}(z) - \ov{\nu}_{k-1}P_{k-1}(z), \quad k\geq1,
$$
and take
$$
\nu_k= \frac{P_k(0)}{Q_k(0)}, \quad k\geq0.
$$ 
Then        
$$
\textup{det} \D_k(P,Q) = (1-|\nu_{k-1}^2|)|Q_{k-1}(0)|^2, \quad k\geq1, 
$$
therefore condition (\ref{STB-eq:SCcriterion.scalar}) holds iff the Schur coefficients of the pair $(P,P^{\#})$ satisfy
\beqq\label{STB-eq:SCcoefcriterion.scalar}
|\nu_k|<1, \quad k=0, \ldots, 3.
\eeqq       
In particular the Schur coefficients read
\beqq\label{STB-eq:SCcoef.scalar}
\nu_0 = p_4, \quad \nu_1=\frac{p_3 - p_4p_1}{1-|p_4|^2}, \quad \nu_2=\frac{(1-|p_4|^2)(p_2 - p_4p_2) - (p_3 - p_4p_1)(p_1-p_4p_3)}{(1-|p_4|^2)^2- |p_3 - p_4p_1|^2},
\eeqq       
\beqq\label{STB-eq:SCcoef4.scalar}
 \nu_3=\frac{p_1-p_4p_3-\frac{(p_3-p_4p_1)(p_2-p_2p_4)}{1-|p_4|^2}-\frac{(1-|p_4|^2)(p_2 - p_4p_2) - (p_3 - p_4p_1)(p_1-p_4p_3)}{(1-|p_4|^2)^2- |p_3 - p_4p_1|^2}(p_1-p_4p_3-\frac{(p_3-p_4p_1)(p_2-p_4p_2)}{1-|p_4|^2})}{1-|p_4|^2-\frac{|p_3-p_4p_1|^2}{1-|p_4|^2}-\frac{(1-|p_4|^2)(p_2 - p_4p_2) - (p_3 - p_4p_1)(p_1-p_4p_3)}{(1-|p_4|^2)^2- |p_3 - p_4p_1|^2}(p_2-p_4p_2-\frac{(p_3-p_4p_1)(p_1-p_4p_3)}{1-|p_4|^2})}
\eeqq       
and thus (\ref{STB-eq:SCcoefcriterion.scalar}) becomes
$$
\left\{ \begin{array}{l}
        |p_4|<1, \\
        |p_3 - p_4p_1|<1-|p_4|^2,\\
        |(1-|p_4|^2)(p_2 - p_4p_2) - (p_3 - p_4p_1)(p_1-p_4p_3)|<(1-|p_4|^2)^2- |p_3 - p_4p_1|^2,\\
        \left|p_1-p_4p_3-\frac{(p_3-p_4p_1)(p_2-p_2p_4)}{1-|p_4|^2}-\frac{(1-|p_4|^2)(p_2 - p_4p_2) - (p_3 - p_4p_1)(p_1-p_4p_3)}{(1-|p_4|^2)^2- |p_3 - p_4p_1|^2}(p_1-p_4p_3-\frac{(p_3-p_4p_1)(p_2-p_4p_2)}{1-|p_4|^2})\right|\\
        <\left|1-|p_4|^2-\frac{|p_3-p_4p_1|^2}{1-|p_4|^2}-\frac{(1-|p_4|^2)(p_2 - p_4p_2) - (p_3 - p_4p_1)(p_1-p_4p_3)}{(1-|p_4|^2)^2- |p_3 - p_4p_1|^2}(p_2-p_4p_2-\frac{(p_3-p_4p_1)(p_1-p_4p_3)}{1-|p_4|^2})\right|
        \end{array}\right.
$$
and simplifying the last condition and using $p_4>0$ we get
\beqq\label{STB-eq:generalSC.scalar}
(SC) \left\{ \begin{array}{l}
        p_4<1, \\
        |p_3 - p_4p_1|<1- (p_4)^2,\\
        |(1-(p_4)^2)(p_2 - p_4p_2) - (p_3 - p_4p_1)(p_1-p_4p_3)|<(1-(p_4)^2)^2- |p_3 - p_4p_1|^2,\\
       |(1+p_4)(p_1- p_4p_3) -p_2(p_3-p_4p_1)|<(1-(p_4)^2)(1+p_2 + p_4) -(p_1+p_3)(p_3 - p_4p_1)        \end{array}\right.
\eeqq
The Schur-Cohn criterion simplifies to  (cf. \cite{jury:1991})
\beqq\label{STB-eq:realSC.scalar}
(SCJ) \left\{ \begin{array}{l}
p_4<1,\\
|p_1+p_3|<1+p_2+p_4,\\
\left|p_2(1-p_4)(1-(p_4)^2) - (p_3-p_4p_1)(p_1-p_4p_3)\right|<(1-(p_4)^2)^2 -(p_3-p_1p_4)^2.
\end{array}\right.
\eeqq
An alternative condition for (SCJ3) reads   \cite[Ex 5.1, p. 255]{elaydi:2005}
$$
\left|p_2(1-p_4) + p_4(1-(p_4)^2) + p_1(p_4p_1-p_3)\right|<p_2p_4(1-p_4) + 1-(p_4)^2 +p_3(p_1p_4-p_3).
$$

\section{Stability conditions for two-step Maruyama methods to the scalar test equation.}\label{NSF:sec:conditions}
\setcounter{equation}{0}

Using the definition of the real coefficients (\ref{STB-eq:charpolcoef1Smat.scalar}) and (\ref{STB-eq:charpolcoef2Smat.scalar}) and the general conditions  (\ref{STB-eq:realSC.scalar})  we can argue when a two-step Maruyama method is asymptotically mean-square stable. In all the following we take $\alpha_0=\gamma_1=1$ and using (\ref{STB-eq:twostepappcoef}) and (\ref{STB-eq:twostepapp1coef}) rewrite the complex coefficients $a, b, c, d$  for the standard schemes
\beqq\label{STB-eq:twostepappcoefnorm}
a=\frac{-\alpha_1 + \beta_1x}{1-\beta_0x}, \quad b=\frac{y}{1-\beta_0x}, \quad c=\frac{-\alpha_2 + \beta_2x}{1-\beta_0x}, \quad d=\frac{\gamma_2y}{1-\beta_0x}
\eeqq
and $b^*, d^*$ for the improved ones
\beqq\label{STB-eq:twostepapp1coefnorm}
b^*= b + \frac{(1+\eta_1)xy}{1-\beta_0x}, \quad d^* = d + \frac{(\gamma_2+\eta_2)xy}{1-\beta_0x},
\eeqq
where also we have used 
$$
x:=h\lam, \qquad     y:=\mu\sqrt{h}.
$$

\bth\label{STB-thr:2step}
The two-step stochastic linear difference equation (\ref{STB-eq:scalar.diff}) is asymptotically mean-square stable iff
\beqq\label{STB-eq:scalar2astep.cond}
|c|^2(|c|^2 + |d|^2)<1,
\eeqq
\beqq\label{STB-eq:scalar2bstep.cond}
|a|^2(1+|c|^2) + |b|^2(1-|c|^2) +2 \Re(\ov{a}bc\ov{d})<(1-|c|^2)^2  - (1-|c|^2)|d|^2- 2\Re(ab\ov{d})- 2\Re(a^2\ov{c})
\eeqq
and
\beam\nonumber
&&\left|\left(-2|c|^2 -|d|^2-2\Re(ab\ov{d}) - 2\Re (a^2\ov{c})\right)(1-|c|^2(|c|^2 + |d|^2))(1-|c|^4(|c|^2 + |d|^2)^2)\right.\\
\nonumber&&\left. \quad- \left(-2\Re(\ov{a}bc\ov{d}) + |c|^2(|b|^2- |a|^2)+|c|^2(|c|^2 + |d|^2)(|a|^2 + |b|^2)\right)\right.\\
\nonumber&&\left.\quad\times(-|a|^2 - |b|^2+|c|^2(|c|^2 + |d|^2)\left(2\Re(\ov{a}bc\ov{d}) - |c|^2(|b|^2- |a|^2)\right)\right|\\
\label{STB-eq:scalar2step2.cond}&&<(1-|c|^4(|c|^2 + |d|^2)^2)^2 -(2\Re(\ov{a}bc\ov{d}) - |c|^2(|b|^2- |a|^2)-|c|^2(|c|^2 + |d|^2)(|a|^2 + |b|^2))^2.
\eeam
Moreover, if  (\ref{STB-eq:scalar2bstep.cond}) holds along with
\beqq\label{STB-eq:scalar2astep2.cond}
|c|^2 + |d|^2<1
\eeqq
and
\beqq\label{STB-eq:scalar2step3.cond}
\Re(ab\ov{d})  + \Re(\ov{a}bc\ov{d})\geq0, \quad \Re(a^2\ov{c})\geq -|a|^2|c|^2,
\eeqq
then (\ref{STB-eq:scalar2step2.cond}) is also true.
For the improved version we take $b^*$ and $d^*$ in place of $b$ and $d$ respectively.
\ethe

\bpf[Proof of Theorem \ref{STB-thr:2step}]
We rewrite the coefficients $p_i, i=1,...,4,$ by (\ref{STB-eq:charpolcoef1Smat.scalar}) and (\ref{STB-eq:charpolcoef2Smat.scalar})   
$$
p_1=-|a|^2 - |b|^2, \quad p_2=-2|c|^2 -|d|^2 -2\Re(ab\ov{d}) - 2\Re (a^2\ov{c}),
$$
$$
p_3=-2\Re(\ov{a}bc\ov{d}) + |c|^2(|b|^2- |a|^2), \quad p_4=|c|^2(|c|^2 + |d|^2).
$$
We need to check conditions (\ref{STB-eq:realSC.scalar}) to conclude about the stability of the method. Condition (SCJ1) implies $|c|^2(|c|^2 + |d|^2)<1.$ Note that
\beam\nonumber
p_1 - p_3 &=&  -|b|^2- |b|^2|c|^2- |a|^2+|a|^2|c|^2 + 2\Re(\ov{a}bc\ov{d})\\
\nonumber&\leq&-|b|^2- |b|^2|c|^2- |a|^2+|a|^2|c|^2 + 2|a||b||c||d|\\
\nonumber&\leq&-|b|^2- |b|^2|c|^2- |a|^2+|a|^2|c|^2 + |a|^2|c|^2|d|^2 + |b|^2\\
\nonumber&<&-|a|^2+|a|^2|c|^2 + |a|^2|d|^2\\
\label{STB-eq:p1-p3.scalar}&=&|a|^2(|c|^2 + |d|^2-1)<0
\eeam
and 
\beam\nonumber
-p_1 - p_3 &=&  |a|^2(1+|c|^2)+|b|^2(1-|c|^2)  + 2\Re(\ov{a}bc\ov{d}) \\
\nonumber&\geq&|a|^2+|a|^2|c|^2 +|b|^2 - |b|^2|c|^2 - 2|a||b||c||d|\\
\nonumber&\geq& |a|^2+|a|^2|c|^2 +|b|^2 - |b|^2|c|^2 - |a|^2 - |b|^2|c|^2|d|^2 \\
\nonumber&>& |b|^2 - |b|^2|c|^2 - |b|^2|d|^2  \\
\label{STB-eq:-p1-p3.scalar}&=&|b|^2(1-|c|^2 - |d|^2)>0,
\eeam
which give 
\beqq\label{STB-eq:p1p3.scalar}
|p_3|<-p_1.
\eeqq
We also have
\beao
|p_1+p_3|&=&\left| -|a|^2 - |b|^2 + |c|^2(|b|^2- |a|^2) -  2\Re(\ov{a}bc\ov{d}) \right|\\
&=&|a|^2(1+|c|^2) + |b|^2(1-|c|)(1+|c|) + 2\Re(\ov{a}bc\ov{d}),
\eeao
due to (\ref{STB-eq:-p1-p3.scalar}) and
\beao
1+p_2+p_4&=&1 -2|c|^2 -|d|^2 -2\Re (a^2\ov{c})   -2\Re(ab\ov{d})  +|c|^4   + |c|^2|d|^2\\
&=&(1-|c|^2)^2  -2\Re (a^2\ov{c})   -2\Re(ab\ov{d})  + (|c|^2-1)|d|^2,
\eeao
so (SCJ2) holds when  (\ref{STB-eq:scalar2bstep.cond}) holds. 
Condition (SCJ3) is (\ref{STB-eq:scalar2step2.cond}).

Furthermore, (\ref{STB-eq:scalar2bstep.cond}), (\ref{STB-eq:scalar2step3.cond}) and (\ref{STB-eq:scalar2astep2.cond}) imply
\beqq\label{STB-eq:for_p1p3p4}
|a|^2 + |b|^2< |a|^2\frac{1+|c|^2}{1-|c|^2} + \frac{2\Re(a^2\ov{c})}{1-|c|^2} + |b|^2<1-|c|^2-|d|^2<1-|c|^2(|c|^2 + |d|^2), 
\eeqq
or
$$
-p_1<1-p_4
$$
which combined with (\ref{STB-eq:p1p3.scalar}) implies
\beqq\label{STB-eq:p1p3p4.scalar}
|p_3-p_4p_1|<1 - (p_4)^2. 
\eeqq
Denote the left-hand side of (\ref{STB-eq:scalar2step2.cond}) by $|L|$ and the right side by $R;$ thus we have to show that $|L|<R.$ Using  (SCJ2) and (\ref{STB-eq:p1p3.scalar}) we get
\beao
L+R&=& (p_2(1-p_4) + 1-(p_4)^2) (1-(p_4)^2) - (p_3-p_4p_1)(p_3-p_4p_1 + p_1-p_4p_3)\\
&=& (p_2 + p_4+1)(1-(p_4)^2)(1-p_4)  -  (p_3-p_4p_1)(p_1+p_3)(1-p_4)\\
&>& -(p_1+p_3)(1-(p_4)^2)(1-p_4)  -  (p_3-p_4p_1)(p_1+p_3)(1-p_4) \\
&=&  -(p_1+p_3)(1-p_4)\left( 1-(p_4)^2 + p_3 - p_4p_1\right)>0,
\eeao 
by (\ref{STB-eq:p1p3p4.scalar}). Therefore $L>-R.$  It remains to prove $L<R.$ 
\beao
L-R&=&p_2(1-p_4)(1-(p_4)^2) - (1-(p_4)^2)^2 - (p_3-p_4p_1)(p_1-p_4p_3) + (p_3-p_1p_4)^2\\
&=&\left(\frac{p_2}{1+p_4}-1\right)(1-(p_4)^2)^2 + (p_3-p_4p_1)(p_3-p_1)(1+p_4)\\
&<&(1-(p_4)^2)(1+p_4)\left( (p_2-p_4-1)\frac{1-p_4}{1+p_4} + (p_3-p_1)\right)\\
&<&(1+p_4)(1-(p_4)^2)\left(p_2-p_4-1+p_3-p_1 -\frac{2p_4}{1+p_4}(p_2-p_4-1) \right)\\
&<&(1+p_4)(1-(p_4)^2)\left(p_2 +3p_4-1+p_3-p_1\right),\\
\eeao
by (\ref{STB-eq:p1-p3.scalar}), (\ref{STB-eq:p1p3p4.scalar})  and (SCJ2). 
The above is negative if  the last term is negative, or equivalently if $p_2+3p_4+p_3-p_1<1.$
 We have
 \beao
&&p_2 + p_3-p_1 +3p_4< -2|c|^2 - |d|^2- 2\Re(a^2\ov{c}) + |c|^2(|b|^2- |a|^2) + |a|^2 + |b|^2 +3 |c|^4 + 3|c|^2|d|^2 \\ 
&<&(|c|^2+1)(|b|^2 + |a|^2) - 2|c|^2  + 3 |c|^4  + 3|c|^2|d|^2 - |d|^2\\
&<& (|c|^2+1)(1-|c|^2-|d|^2) + 2|c|^2(|c|^2 +|d|^2 -1)+ |c|^4 - (1-|c|^2)|d|^2<1,
\eeao
where we used (\ref{STB-eq:scalar2astep2.cond}), (\ref{STB-eq:scalar2step3.cond}) and  (\ref{STB-eq:for_p1p3p4}).
\epf

\subsection{Two-step Adams-Bashforth and Adams-Moulton Maruyama methods.}

In this case $\gamma_2=0$ and thus $d=0$ so the recurrences (\ref{STB-eq:scalar.diff}) simplify to 
\beqq\label{STB-eq:scalarABAM.diff}
X_{i} = aX_{i-1} + cX_{i-2} + bX_{i-1}\xi_{i-1}, \, i=2,3,\ldots,
\eeqq
for the standard schemes and to 
\beqq\label{STB-eq:scalarABIAMI.diff}
X_{i} = aX_{i-1} + cX_{i-2} + b^*X_{i-1}\xi_{i-1}  + d^*X_{i-2}\xi_{i-2}, \, i=2,3,\ldots,
\eeqq
for the improved ones.

\bpr\label{STB-prop:ABAM}
The two-step stochastic linear difference equation (\ref{STB-eq:scalarABAM.diff}) is asymptotically mean-square stable iff
\beqq\label{STB-eq:scalarABAM.cond}
|c|<1, \quad  |a|^2(1+|c|^2) + 2\Re(a^2\ov{c})<(1-|c|^2)^2, 
\eeqq
\beqq\label{STB-eq:scalar2ABAM.cond}
|b|^2<1-|c|^2 - |a|^2\frac{1+|c|^2}{1-|c|^2} - 2\frac{\Re(a^2\ov{c})}{1-|c|^2},
\eeqq
and
\beqq\label{STB-eq:scalar3ABAM.cond}
\Re(a^2\ov{c})\geq -|a|^2|c|^2,
\eeqq
whereas the two-step stochastic linear difference equation (\ref{STB-eq:scalarABIAMI.diff}) is  asymptotically mean-square stable iff  conditions (\ref{STB-eq:scalar2astep.cond}),(\ref{STB-eq:scalar2bstep.cond}) and  (\ref{STB-eq:scalar2step2.cond})  hold or conditions (\ref{STB-eq:scalar2astep2.cond}),(\ref{STB-eq:scalar2bstep.cond}) and  (\ref{STB-eq:scalar2step3.cond}) hold where $b^*$ and $d^*$ are replacing $b$ and $d$ respectively.
\epr 

\bpf[Proof of Proposition \ref{STB-prop:ABAM}]
We show the first case since the second one is a direct application of Theorem \ref{STB-thr:2step}. In the case of   (\ref{STB-eq:scalarABAM.diff}) the coefficients  read  
\beqq\label{STB-eq:charpolcoef1SmatABAM.scalar}
p_1=-|a|^2 - |b|^2, \quad p_2=-2|c|^2 -2\Re (a^2\ov{c}),
\eeqq
\beqq\label{STB-eq:charpolcoef2SmatABAM.scalar}
p_3=|c|^2(|b|^2- |a|^2), \quad p_4=|c|^4.
\eeqq
We apply Theorem \ref{STB-thr:2step} when $d=0.$  
Conditions (\ref{STB-eq:scalar2astep.cond}) or (\ref{STB-eq:scalar2astep2.cond}) are equivalent to $|c|<1.$ Condition (\ref{STB-eq:scalar2bstep.cond}) is just the right-side of (\ref{STB-eq:scalarABAM.cond}) and (\ref{STB-eq:scalar2ABAM.cond}). Finally  (\ref{STB-eq:scalar2step3.cond}) shrinks to $\Re(a^2\ov{c})\geq -|a|^2|c|^2.$
\epf

\bre\label{STB-rem:ABAM}
Consider the case $a, b, c\in\bbR.$ Then conditions (\ref{STB-eq:scalarABAM.cond}), (\ref{STB-eq:scalar2ABAM.cond}) and (\ref{STB-eq:scalar3ABAM.cond}) read (see also \cite[Cor. 6]{tocino_senosiain:2014})
\beqq\label{STB-eq:scalar4ABAM.cond}
0<c<1, \quad |a|<1-c, 
\eeqq
\beqq\label{STB-eq:scalar5ABAM.cond}
\quad b^2(1-c)<(1+c)\left((1-c)^2-a^2\right). 
\eeqq
\ere

\subsection{Schemes for hereditary systems.}

Hereditary systems are used to model processes in a variety of fields such as physics, biology, economy, just to name a few, (cf. \cite{kolmanovski_myshkis:92}). Due to their applications, we present them in a separate subsection. The following stochastic difference equation was  proposed in \cite{shaikhet:97}, 
\beqq
X_{i+1}=\sum_{j=0}^k\alpha_jX_{i-j} + \sigma X_{i-l}\xi_i,
\eeqq
where necessary and sufficient conditions were given  concerning their  asymptotic mean-square stability of the zero solution. By taking the trivial case $l=0$ of this delay system with $k=2$  this falls in our setting (\ref{STB-eq:scalar.diff})  with $b=0,$ that is,
\beqq\label{STB-eq:scalarHER.diff}
X_{i} = aX_{i-1} + cX_{i-2} + dX_{i-2}\xi_{i-2}, \, i=2,3,\ldots,
\eeqq
for the standard schemes and to 
\beqq\label{STB-eq:scalarHERI.diff}
X_{i} = aX_{i-1} + cX_{i-2} + b^*X_{i-1}\xi_{i-1}  + d^*X_{i-2}\xi_{i-2}, \, i=2,3,\ldots,
\eeqq
for the improved ones. 

\bpr\label{STB-prop:HER}
The two-step stochastic linear difference equation (\ref{STB-eq:scalarHER.diff}) is asymptotically mean-square stable iff
\beqq\label{STB-eq:scalarHER.cond}
|c|^2  + |d|^2<1, \quad |a|^2(1+|c|^2) +2\Re(a^2\ov{c})<(1-|c|^2)^2  - (1-|c|^2)|d|^2 
\eeqq
and
\beqq\label{STB-eq:scalar2HER.cond}
\Re(a^2\ov{c})\geq -|a|^2|c|^2,\eeqq
whereas the two-step stochastic linear difference equation (\ref{STB-eq:scalarHERI.diff}) is  asymptotically mean-square stable iff  conditions (\ref{STB-eq:scalar2astep.cond}),(\ref{STB-eq:scalar2bstep.cond}) and  (\ref{STB-eq:scalar2step2.cond})  hold or conditions (\ref{STB-eq:scalar2astep2.cond}),(\ref{STB-eq:scalar2bstep.cond}) and  (\ref{STB-eq:scalar2step3.cond}) hold where $b^*$ and $d^*$ are replacing $b$ and $d$ respectively.
\epr 

\bpf[Proof of Proposition \ref{STB-prop:HER}]
We show the first case since the second one is a direct application of Theorem \ref{STB-thr:2step}. In the case of   (\ref{STB-eq:scalarHER.diff}) the coefficients  read  
\beqq\label{STB-eq:charpolcoef1SmatHER.scalar}
p_1=-|a|^2, \quad p_2=-2|c|^2 -|d|^2-2\Re (a^2\ov{c}),
\eeqq
\beqq\label{STB-eq:charpolcoef2SmatHER.scalar}
p_3=- |a|^2|c|^2, \quad p_4=|c|^2(|c|^2 + |d|^2).
\eeqq
We apply Theorem \ref{STB-thr:2step} when $b=0.$  
Conditions (\ref{STB-eq:scalar2astep.cond}) or (\ref{STB-eq:scalar2astep2.cond}) are equivalent to $|c|^2+|d|^2<1.$ Condition (\ref{STB-eq:scalar2bstep.cond}) is just the right-side of (\ref{STB-eq:scalarHER.cond}). Finally  (\ref{STB-eq:scalar2step3.cond}) shrinks to $\Re(a^2\ov{c})\geq -|a|^2|c|^2.$
\epf

\bre\label{STB-rem:HER}
Consider the case $a, c, d\in\bbR.$ Then conditions (\ref{STB-eq:scalarHER.cond}) and (\ref{STB-eq:scalar2HER.cond})  read (see also \cite[Cor. 5]{tocino_senosiain:2014})
\beqq\label{STB-eq:scalar4HER.cond}
|c|^2 + |d|^2<1, \quad |a|<1-c, 
\eeqq
\beqq\label{STB-eq:scalar5HER.cond}
\quad \frac{1-c}{(1+c)((1-c)^2-a^2)}<\frac{1}{d^2}. 
\eeqq
\ere

\section{Linear MS-stability.}\label{NSF:sec:linearMSstability}
\setcounter{equation}{0}

Recall the scalar linear test-equation  (\ref{STB-eq:sde.scalar})

$$
dX(t) = \lam X(t) dt + \mu X(t) dW(t), \quad X(t_0)=X_0,
$$
where $\lam, \mu, X_0\in\bbC.$ Its  zero solution is asymptotically mean-square stable iff $\Re(\lam) + |\mu|^2/2<0;$ in the case $\mu=0$ the above condition reduces to the notion of A-stability.
 The set 
$$
\bbs_{SDE}=\{(\lam, \mu)\in\bbC\times\bbC:  \Re(\lam) + \frac{|\mu|^2}{2}<0\},
$$
is called the mean-square (MS-)stability domain of the stochastic equation (\ref{STB-eq:sde.scalar}). In an analogous manner the MS-stability domain of a two-step stochastic method (SM) for a given step size $h>0$ is defined as  
\beqq\label{STB-eq:sde.domain_method}
\bbs_{SM}(h)=\{(\lam, \mu)\in\bbC\times\bbC: \hbox{conditions } (\ref{STB-eq:scalar2astep.cond}), (\ref{STB-eq:scalar2bstep.cond}) \hbox{ and } (\ref{STB-eq:scalar2step2.cond}) \hbox{ hold}\}.
\eeqq

In case $\lam, \mu\in\bbR$ we have the notions of the stability regions 
\beqq\label{STB-eq:sde.region_sde}
\bbr_{SDE}=\{(\lam, \mu)\in\bbR\times\bbR:  \lam  + \frac{\mu^2}{2}<0\},
\eeqq
for the sde and 
\beqq\label{STB-eq:sde.region_method}
\bbr_{SM}(h)=\{(\lam, \mu)\in\bbR\times\bbR: \hbox{conditions } (\ref{STB-eq:scalar2astep.cond}), (\ref{STB-eq:scalar2bstep.cond})  \hbox{ and } (\ref{STB-eq:scalar2step2.cond}) \hbox{ hold}\}
\eeqq
for the method. A stochastic method is said to be MS-stable if 
$$
\bbs_{SDE}\subseteq \bbs_{SM}, \mbox{ or } \bbr_{SDE}\subseteq\bbr_{SM}   \mbox{ for all } h>0.
$$
The inverse relation
$$
\bbs_{SM}\subset \bbs_{SDE}, \mbox{ or } \bbr_{SM}\subset\bbr_{SDE}   \mbox{ for all } h>0.
$$
means that the method is unstable whenever the test-equation is unstable. In this case the notion of conditional MS-stability comes to play, where one has to determine a step size $h_0$ such that for a given pair of $(\lam, \mu)$ in the stability domain or region of the sde the method is mean-square stable for all $h<h_0.$

\subsection{MS-stability of Adams-Bashforth Maruyama scheme.}
The coefficients of the AB2 scheme, see Table \ref{STB-tab:parameters}, read
$$
a=1+\frac{3}{2}x, \quad b=y, \quad c=-\frac{1}{2}x, \quad d=0
$$
and for the improved AB2I
$$
b^*= y(1+x), \quad d^* = -\frac{1}{2}xy.
$$

\begin{table}\caption{Parameters of two-step schemes as in (\ref{STB-eq:scalar.diff}).}\label{STB-tab:parameters}
 \begin{tabular}{r|| r r r r} 
  \hline
Method & $a$ & $b$ & $c$ & $d$\\ [0.5ex] 
 \hline\hline
 AB2 & $1+(3/2)x$ & $y$ & $-x/2$ & $0$\\ [0.5ex]
 \hline
 AB2I & $1+(3/2)x$ & $y(1+x)$ & $-x/2$ & $-xy/2$\\[0.5ex]
 \hline
 AM2 & $\frac{1+(8/12)x}{1 - (5/12)x}$ & $\frac{y}{1 - (5/12)x}$ & $\frac{-x/12}{1 - (5/12)x}$  & 0\\[0.5ex]
 \hline
 AM2I & $\frac{1+(8/12)x}{1 - (5/12)x}$ &$\frac{y+(7/12)xy}{1 - (5/12)x}$ & $\frac{-x/12}{1 - (5/12)x}$ & $\frac{-xy/12}{1 - (5/12)x}$\\[0.5ex]
 \hline
BDF2 & $\frac{4/3}{1 - (2/3)x}$ & $\frac{y}{1 - (2/3)x}$ & $\frac{-1/3}{1 - (2/3)x}$ & $\frac{-y/3}{1 - (2/3)x}$\\ [0.5ex]
 \hline
BDF2I & $\frac{4/3}{1 - (2/3)x}$ & $\frac{y+ xy/3}{1 - (2/3)x}$ &$\frac{-1/3}{1 - (2/3)x}$ & $\frac{-y/3}{1 - (2/3)x}$\\  [0.5ex]
\end{tabular}
\end{table}

First we take $(\lam,\mu)\in \bbs_{AB2}$ where
\beqq\label{STB-eq:sde.domain_AB2}
\bbs_{AB2}(h)=\{(\lam, \mu)\in\bbC\times\bbC: \hbox{conditions } (\ref{STB-eq:scalarABAM.cond}), (\ref{STB-eq:scalar2ABAM.cond}) \hbox{ and } (\ref{STB-eq:scalar3ABAM.cond}) \hbox{ hold}\}.
\eeqq
Conditions (\ref{STB-eq:scalarABAM.cond}) give
\beqq\label{STB-eq:sde.cond_AB2}
|x|<2, \quad  \left(1+\frac{|x|^2}{4}\right)\left|1+\frac{3}{2}x\right|^2 -\Re\left((1+\frac{3}{2}x)^2\ov{x}\right)<\left(1-\frac{1}{4}|x|^2\right)^2.
\eeqq
Now, inspecting the second inequality further we conclude that 
\beao
&&\left(1+\frac{|x|^2}{4}\right)\left(1+3\Re(x) + \frac{9}{4}|x|^2\right)  -\Re(\ov{x})-\frac{9}{4}|x|^2\Re(x)-3|x|^2\\
&=&1 -\frac{1}{2}|x|^2 + 2\Re(x) -\frac{3}{2}|x|^2\Re(x) + \frac{9|x|^4}{16}\\
&<& 1 -\frac{1}{2}|x|^2 + \frac{|x|^4}{16},
\eeao
when 
$$
(2 -\frac{3}{2}|x|^2)\Re(x)<-\frac{|x|^4}{2},
$$
which implies $\Re(x)<0,$ that is $\Re(\lam)<0,$ when $|x|^2<\frac{4}{3}.$  Conditions (\ref{STB-eq:scalar2ABAM.cond}) give
$$
|y|^2<\frac{4}{4-|x|^2} \left( -2\Re(x) +\frac{3}{2}|x|^2\Re(x)- \frac{|x|^4}{2}\right)<-2\Re(x),
$$
when
$$
\left( \frac{4-3|x|^2}{4-|x|^2}-1\right)\Re(x) + \frac{|x|^4}{4-|x|^2}>0,
$$
which holds for any $0<|x|<2$ with $\Re(x)<0.$ Moreover, condition (\ref{STB-eq:scalar3ABAM.cond}) reads, 
$$
\Re\left(   -(1+\frac{3}{2}x)^2\frac{\ov{x}}{2}\right)=\frac{1}{2}\left(-(1+\frac{9}{4}|x|^2)\Re(x) - 3|x|^2\right)\geq-\left(1+3\Re(x) + \frac{9}{4}|x|^2\right)\frac{|x|^2}{4},
$$
or equivalently  
$$
 (6|x|^2+8)\Re(x)\leq 9|x|^4 -20|x|^2,
$$
which implies $\Re(x)<0$ when $|x|^2<20/9.$ Conditions (\ref{STB-eq:scalarABAM.cond}), (\ref{STB-eq:scalar2ABAM.cond}) and (\ref{STB-eq:scalar3ABAM.cond}) 
 hold when 
$$
|x|<1, \qquad |y|^2<\frac{2}{4-|x|^2} \left( -4\Re(x) +3|x|^2\Re(x)- |x|^4\right).
$$
Therefore we  get 
$$
\bbs_{AB2}(h)\subset \bbs_{SDE},$$
for any $h>0,$  which means that AB2 is unstable whenever the test-equation is unstable. Now, given 
  $(\lam,\mu)\in\bbC\times\bbC$  we want to find $h_0>0$ such that $ \bbs_{SDE} \subset \bbs_{AB2}(h)$ for any $h<h_0.$ Since we chose the parameters in the stability domain  $\bbs_{SDE}$ we have that $|\mu|^2<-2\Re(\lam).$   The relation $|x|<1$ gives 
  $$h<\frac{1}{|\lam|}.$$ 
Moreover, by (\ref{STB-eq:scalar2ABAM.cond}) we need to show that 
\beao
&&h|\mu|^2 + \frac{2}{4-h^2|\lam|^2} \left(  -3h^3|\lam|^2\Re(\lam) + h^4|\lam|^4
+ 4h\Re(\lam)\right)\\
&=&\frac{h}{4-h^2|\lam|^2}\left(4|\mu|^2-h^2|\lam|^2|\mu|^2 - 6h^2|\lam|^2\Re(\lam) + 2h^3|\lam|^4 + 8\Re(\lam)\right)<0,
\eeao
 which holds when  
 $$
\un{- 6h^2|\lam|^2\Re(\lam) + 4|\mu|^2+ 8\Re(\lam)}_{negative}  +  \un{2h^3|\lam|^4-h^2|\lam|^2|\mu|^2}_{negative}<0,
$$
or in terms of $h$ for 
$$h<\min\left\{ \frac{|\mu|^2}{2|\lam|^2},    \sqrt{\frac{4(-2\Re(\lam)-|\mu|^2)}{-6\Re(\lam)|\lam|^2}}\right\} :=h_1.$$ 
So  given $(\lam,\mu)\in\bbs_{SDE}$   the method AB2 is conditionally MS-stable for any $h<h_0$ where 
$$
h_0=\min \left\{\frac{1}{|\lam|}, h_1\right\}.
$$

In case the parameters are real conditions (\ref{STB-eq:scalarABAM.cond}), (\ref{STB-eq:scalar2ABAM.cond}) and (\ref{STB-eq:scalar3ABAM.cond}) shrink to (\ref{STB-eq:scalar4ABAM.cond}) and (\ref{STB-eq:scalar5ABAM.cond}) respectively by Remark \ref{STB-rem:ABAM}. The asymptotic region reads
$$
\bbr_{AB2}(h)=\left\{(\lam,\mu)\in\bbR\times\bbR: -1<\lam h<0, \mu^2<\frac{2\lam (\lam h-2)(\lam h+1)}{\lam h+2}\right\}
$$
and $\bbr_{AB2}(h)\subset \bbr_{SDE}$ for any $h>0.$  Given $(\lam,\mu)\in\bbr_{SDE}$ and  $h>0$ the method AB2 is conditionally MS-stable for all $h<h_0$ where 
$$
h_0=\min \left\{-\frac{1}{\lam},\frac{\mu^2 + 2\lam + \sqrt{(\mu^2 + 2\lam)(\mu^2 + 18\lam)}}{4\lam^2}\right\}.
$$
 
\begin{figure}
  \centering
  \begin{subfigure}{.45\textwidth} 
   \includegraphics[width=0.8\textwidth]{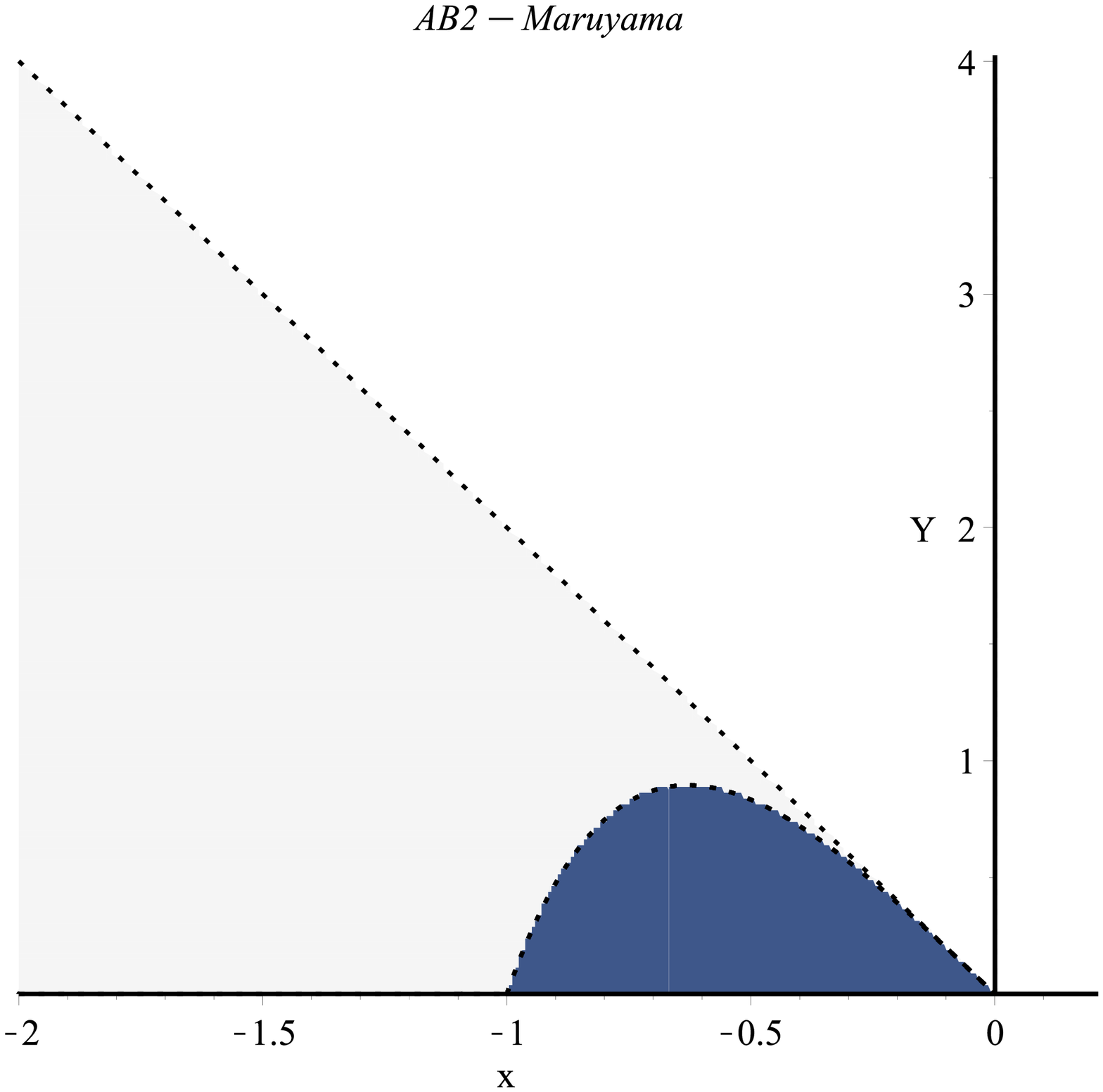}\label{STB-fig:ab2}
     \caption{Stability region of AB2 (Niagara-Blue).}
\end{subfigure}
\begin{subfigure}{.45\textwidth}
   \includegraphics[width=0.8\textwidth]{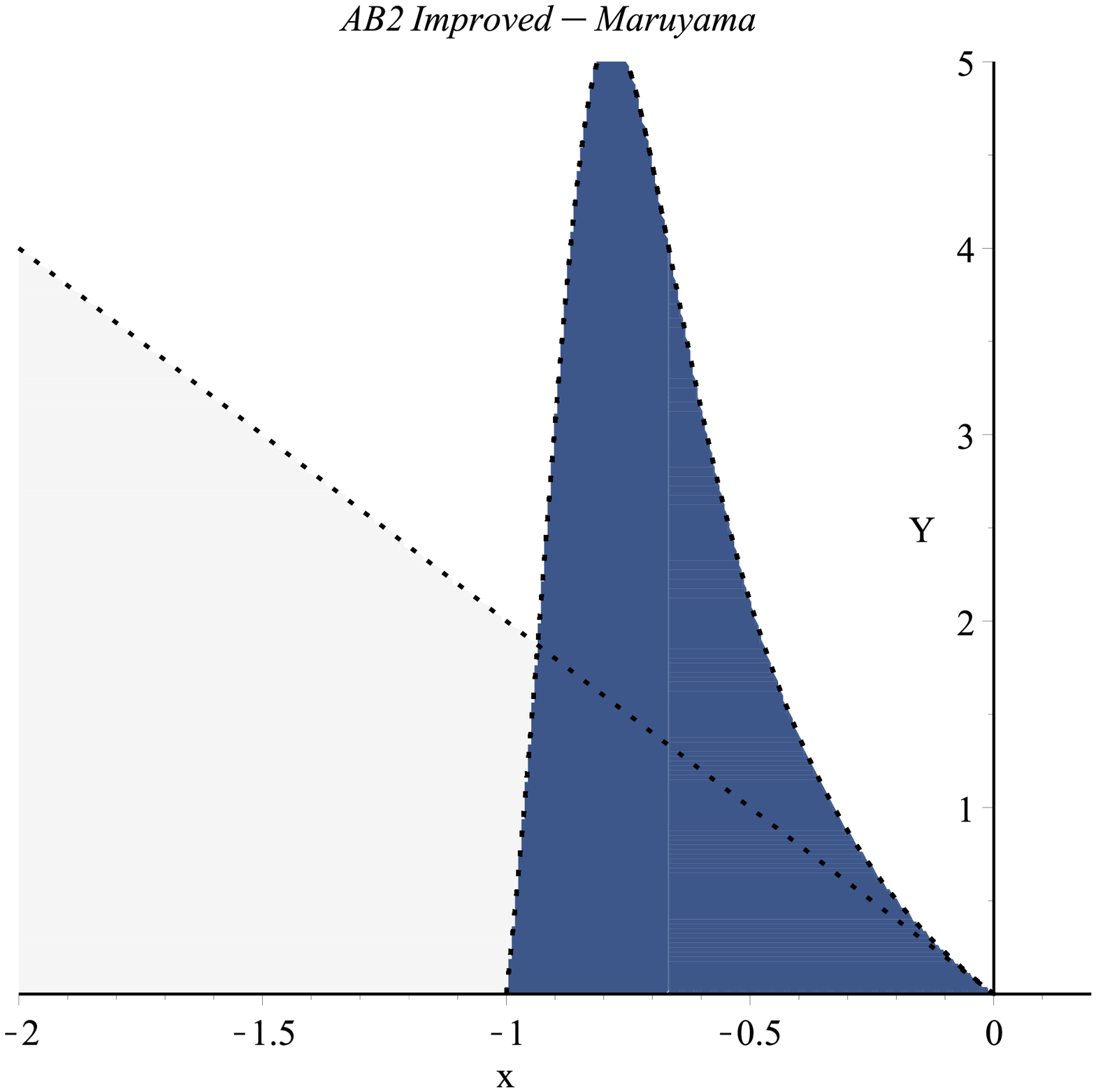}\label{STB-fig:ab2i}
    \caption{Stability region of AB2I (Niagara-Blue).}
\end{subfigure}
\caption{Stability regions of AB2- and Improved AB2-Maruyama.}\label{STB-fig:ab2ab2i} 
  \end{figure}

In Figure \ref{STB-fig:ab2ab2i} we represent the stability regions of the AB2 and AB2I scheme respectively in the $(x,Y)$-plane where $x=\lam h$ and $Y=\mu^2 h,$ where also the stability region of the SDE is shown (it corresponds to the region $0<Y<-2x,$ that is the light-shaded triangle.)
 
\subsection{MS-stability of Adams-Moulton Maruyama scheme.}

The coefficients of the AM2 scheme, see Table \ref{STB-tab:parameters}, read
$$
a=\frac{1+(8/12)x}{1 - (5/12)x}, \quad b=\frac{y}{1 - (5/12)x}, \quad c=-\frac{x/12}{1 - (5/12)x}, \quad d=0
$$
and for the improved AM2I
$$
b^*= \frac{y+(7/12)xy}{1 - (5/12)x}, \quad d^* = -\frac{xy/12}{1 - (5/12)x}.
$$
 First we take $(\lam,\mu)\in \bbs_{AM2}$ where
\beqq\label{STB-eq:sde.domain_AM2}
\bbs_{AM2}(h)=\{(\lam, \mu)\in\bbC\times\bbC: \hbox{conditions } (\ref{STB-eq:scalarABAM.cond}), (\ref{STB-eq:scalar2ABAM.cond}) \hbox{ and } (\ref{STB-eq:scalar3ABAM.cond}) \hbox{ hold}\}.
\eeqq
Conditions (\ref{STB-eq:scalarABAM.cond}) give
\beqq\label{STB-eq:sde.cond_AM2}
|x|<|12-5x|, \,\left(1+\frac{|x|^2}{|12-5x|^2}\right)\frac{|12+8x|^2}{|12 - 5x|^2} -2\Re\left(  \frac{(12+8x)^2}{(12 - 5x)^2}\frac{\ov{x}}{12-5\ov{x}}\right)<\left(1-\frac{|x|^2}{|12-5x|^2}\right)^2.
\eeqq
The first inequality is satisfied by those $x$ with $5\Re(x)<|x|^2+6$
and the second inequality holds when $-6<\Re(x)<0.$ Therefore (\ref{STB-eq:sde.cond_AM2}) holds iff
\beqq\label{STB-eq:sde.condeq_AM2}
-6<\Re(x)<0,
\eeqq
which imply $\Re(\lam)<0.$  Conditions (\ref{STB-eq:scalar2ABAM.cond}) give
$$
|y|^2<\frac{1}{|12-5x|^2-|x|^2} \left( \frac{(|12-5x|^2-|x|^2)^2 - (|12-5x|^2 +|x|^2)|12+8x|^2 }{|12-5x|^2} + 2\Re\left(  \frac{(12+8x)^2}{(12 - 5x)}\ov{x}\right)\right),
$$
which is smaller than $-2\Re(x)$ for any $\Re(x)<0.$ Moreover, condition (\ref{STB-eq:scalar3ABAM.cond}) reads, 
$$
\Re\left( -\frac{(12+8x)^2}{(12-5x)^2}\frac{\ov{x}}{12-5\ov{x}}\right)=-\frac{1}{|12-5x|^4}\Re\left((12+8x)^2\ov{x}(12-5\ov{x})\right)\geq  -\frac{|12 + 8x|^2|x|^2}{|12-5x|^4}
$$
or equivalently
$$
-12\cdot(144-32|x|^2)\Re(x) +384|x|^4 - 15\cdot144|x|^2 + 5\cdot144\Re(\ov{x}^2)\geq0,$$
which implies $\Re(x)<0$ and $|x|^2<9/2.$ Therefore we  get 
$$
\bbs_{AM2}(h)\subset \bbs_{SDE},$$
for any $h>0,$  which means that AM2 is unstable whenever the test-equation is unstable. Now, given   $(\lam,\mu)\in\bbC\times\bbC$  we want to find $h_0>0$ such that $ \bbs_{SDE} \subset \bbs_{AM2}(h)$ for any $h<h_0.$ Since we chose the parameters in the stability domain  $\bbs_{SDE}$ we have that $|\mu|^2<-2\Re(\lam).$   Relation (\ref{STB-eq:sde.condeq_AM2}) implies
 $$h<-\frac{6}{\Re(\lam)}.$$   Moreover, by (\ref{STB-eq:scalar2ABAM.cond}) we need to show that 
$$h|\mu|^2 +  \frac{h^2|\lam|^2}{|12-5h\lam|^2} + \frac{|12-5h\lam|^2   + h^2|\lam|^2}{|12-5h\lam|^2-h^2|\lam|^2}\frac{|12+8h\lam|^2 }{|12-5h\lam|^2} - \frac{2h}{|12-5h\lam|^2-h^2|\lam|^2}\Re\left(  \frac{(12+8h\lam)^2}{(12 - 5h\lam)}\ov{\lam}\right)<1\\
$$
which holds for sufficiently small $h_1>0$  implying that  the method AM2 is conditionally MS-stable for any $h<h_0$ where 
$$
h_0=\min \left\{-\frac{6}{\Re(\lam)}, h_1\right\}.
$$

In case the parameters are real we need to show  (\ref{STB-eq:scalar4ABAM.cond}) and (\ref{STB-eq:scalar5ABAM.cond}) respectively. The asymptotic region reads
$$
\bbr_{AM2}(h)=\left\{(\lam,\mu)\in\bbR\times\bbR: -6<\lam h<0, \mu^2<\frac{\lam (\lam h-2)(\lam h+6)}{2(3 - \lam h)}\right\}
$$
and $\bbr_{AM2}(h)\subset \bbr_{SDE}$ for any $h>0.$  Given $(\lam,\mu)\in\bbr_{SDE}$ and  $h>0$ the method AM2 is conditionally MS-stable for all $h<h_0$ where 
$$
h_0=\min \left\{-\frac{6}{\lam},\frac{-\mu^2 - 2\lam + \sqrt{(\mu^2 + 2\lam)(\mu^2 + 8\lam)}}{\lam^2}\right\}.
$$
 
\begin{figure}
  \centering
  \begin{subfigure}{.45\textwidth} 
   \includegraphics[width=0.8\textwidth]{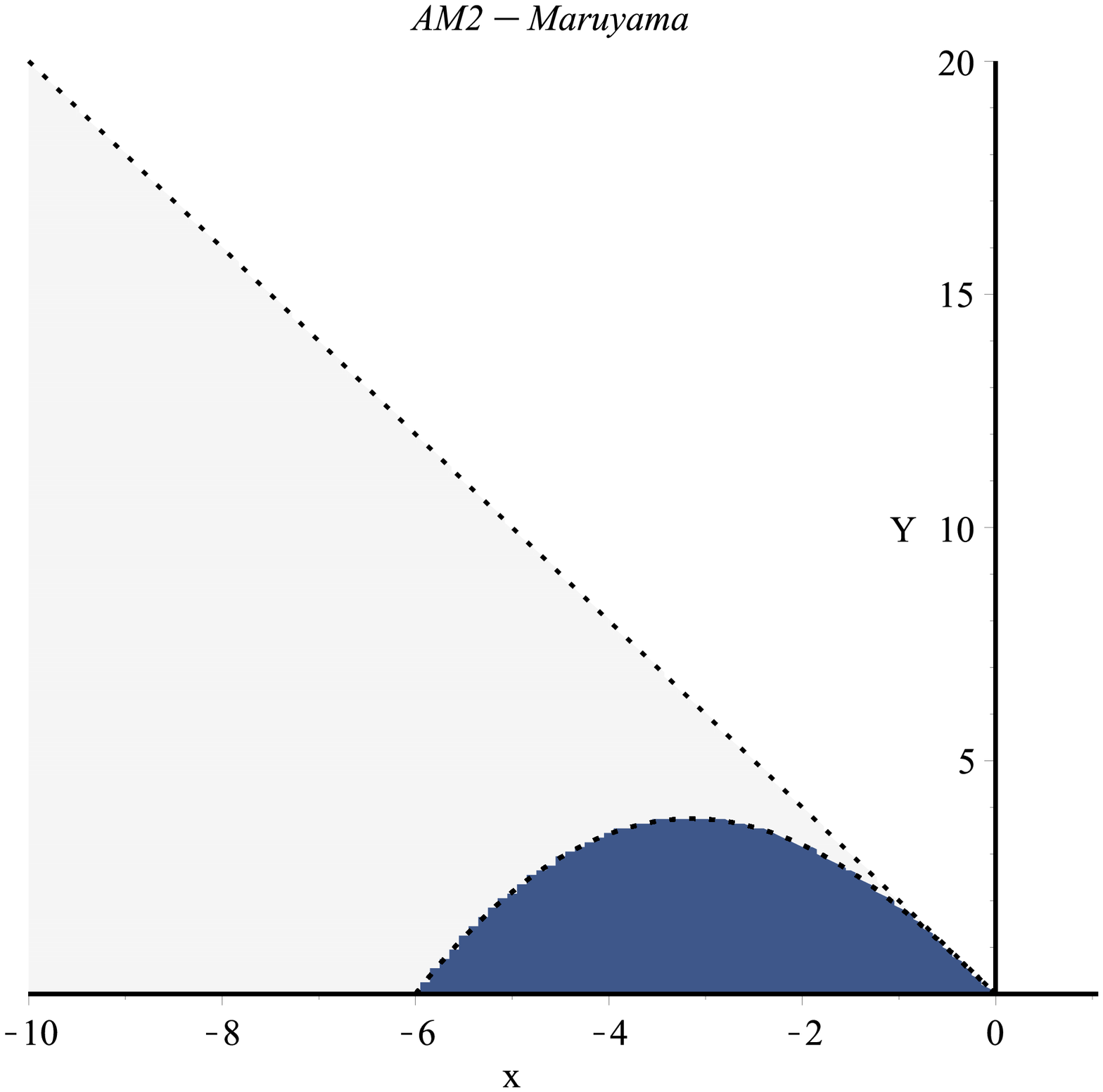}\label{STB-fig:am2}
     \caption{Stability region of AM2 (Niagara-Blue).}
\end{subfigure}
\begin{subfigure}{.45\textwidth}
   \includegraphics[width=0.8\textwidth]{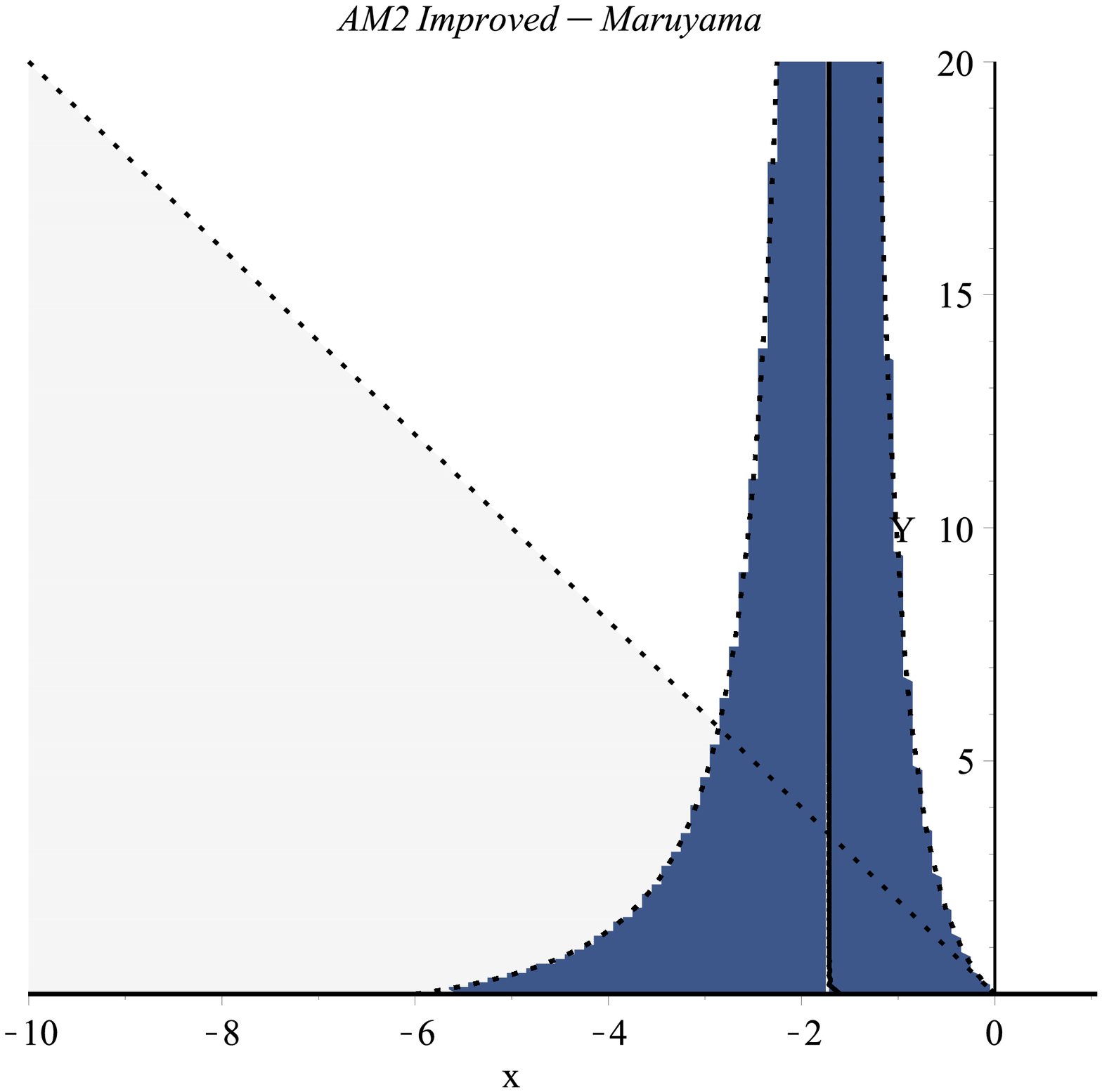}\label{STB-fig:am2i}
    \caption{Stability region of AM2I (Niagara-Blue).}
\end{subfigure}
\caption{Stability regions of AM2- and Improved AM2-Maruyama.}\label{STB-fig:am2am2i} 
  \end{figure}

In Figure \ref{STB-fig:am2am2i} we represent the stability regions of the AM2 and AM2I scheme respectively in the $(x,Y)$-plane where $x=\lam h$ and $Y=\mu^2 h,$ where also the stability region of the SDE is shown (it corresponds to the region $0<Y<-2x,$ that is the light-shaded triangle.)

\subsection{Two-step BDF Maruyama scheme.}
The coefficients of the BDF2 scheme, see Table \ref{STB-tab:parameters}, read
$$
a=\frac{4/3}{1 - (2/3)x}, \quad b=\frac{y}{1 - (2/3)x}, \quad c=-\frac{1/3}{1 - (2/3)x}, \quad d=-\frac{y/3}{1 - (2/3)x}
$$
and for the improved BDF2I
$$
b^*= \frac{y+ xy/3}{1 - (2/3)x}, \quad d^* = -\frac{y/3}{1 - (2/3)x}=d.
$$
 
\begin{figure}
  \centering
  \begin{subfigure}{.45\textwidth} 
   \includegraphics[width=0.8\textwidth]{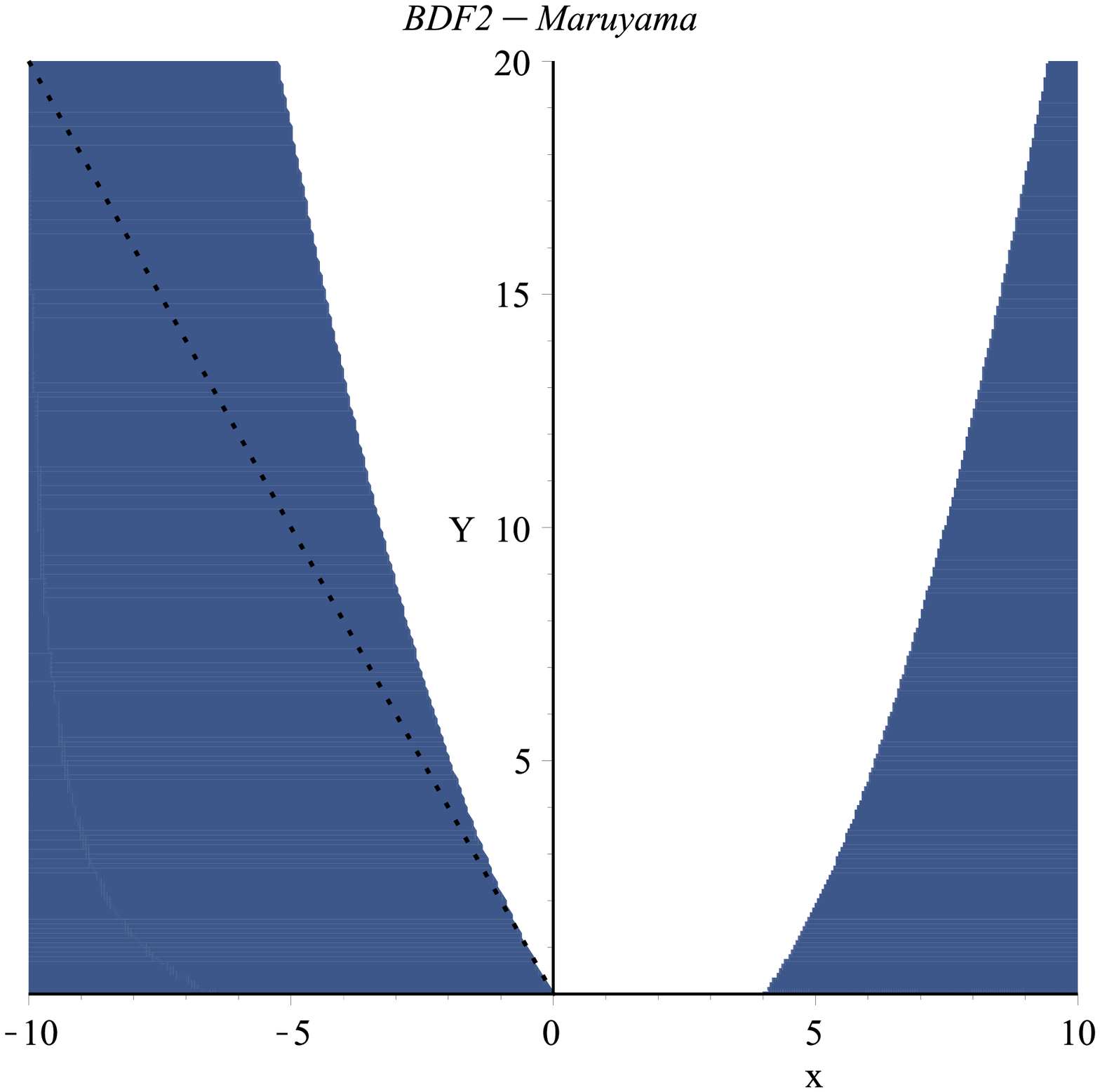}\label{STB-fig:bdf2}
     \caption{Stability region of BDF2 (Niagara-Blue)}
\end{subfigure}
\begin{subfigure}{.45\textwidth}
   \includegraphics[width=0.8\textwidth]{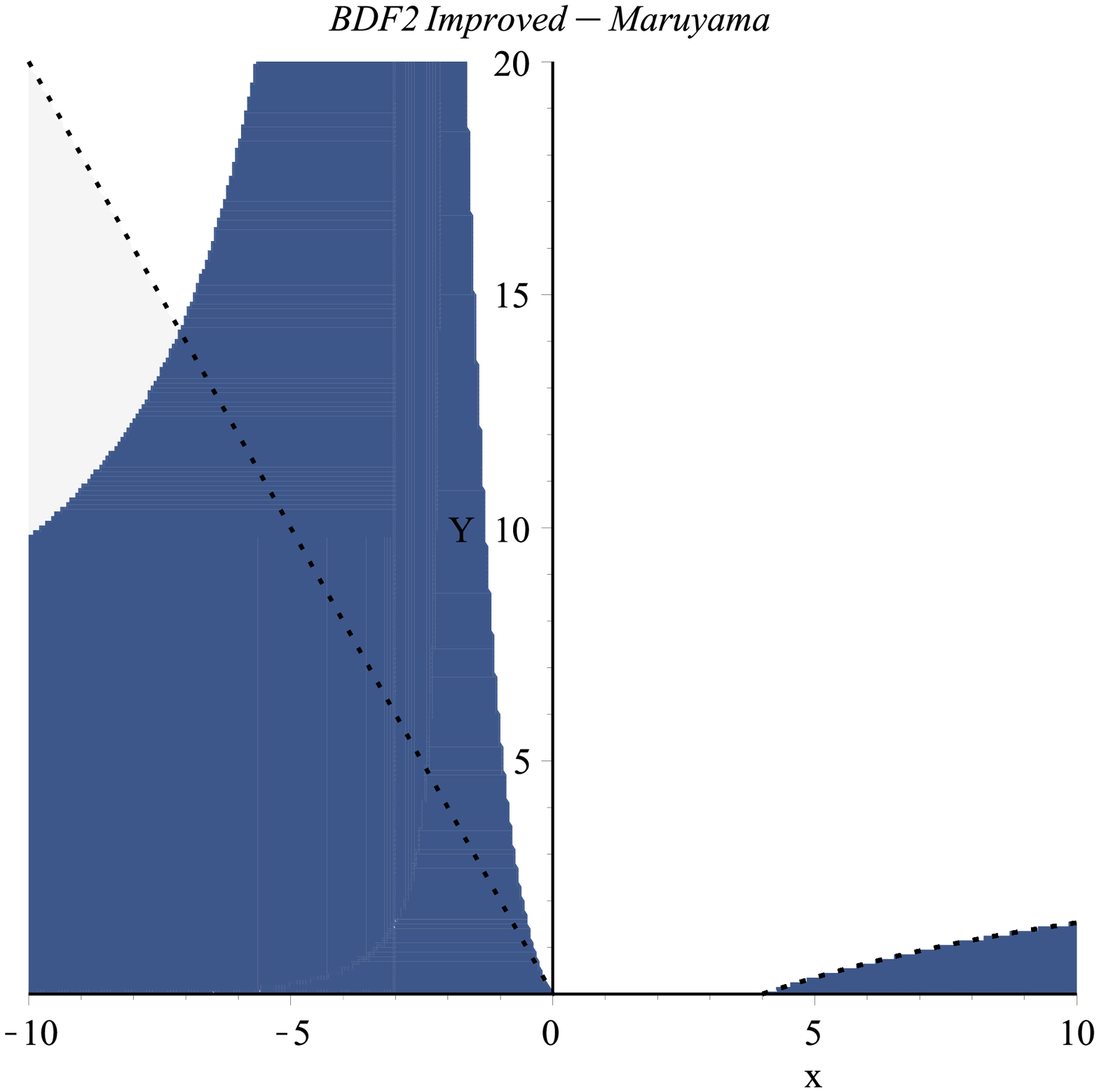}\label{STB-fig:bdf2i}
    \caption{Stability region of BDF2I (Niagara-Blue)}
\end{subfigure}
\caption{Stability regions of BDF2- and Improved BDF2-Maruyama.}\label{STB-fig:bdf2bdf2i} 
  \end{figure}

In Figure \ref{STB-fig:bdf2bdf2i} we represent the stability regions of the BDF2 and BDF2I scheme respectively in the $(x,Y)$-plane where $x=\lam h$ and $Y=\mu^2 h,$ where also the stability region of the SDE is shown (it corresponds to the region $0<Y<-2x,$ that is the light-shaded triangle.)

\section{Experiments.}\label{NSF:sec:exp}
\setcounter{equation}{0}

\begin{figure}[ht]
  \caption{Approximations of the $2^{nd}$ moment of the linear scalar test equation (\ref{STB-eq:sde.scalar}) in the interval $[0,1]$ for different two-step numerical  methods.}
  \centering
   \includegraphics[width=0.8\textwidth]{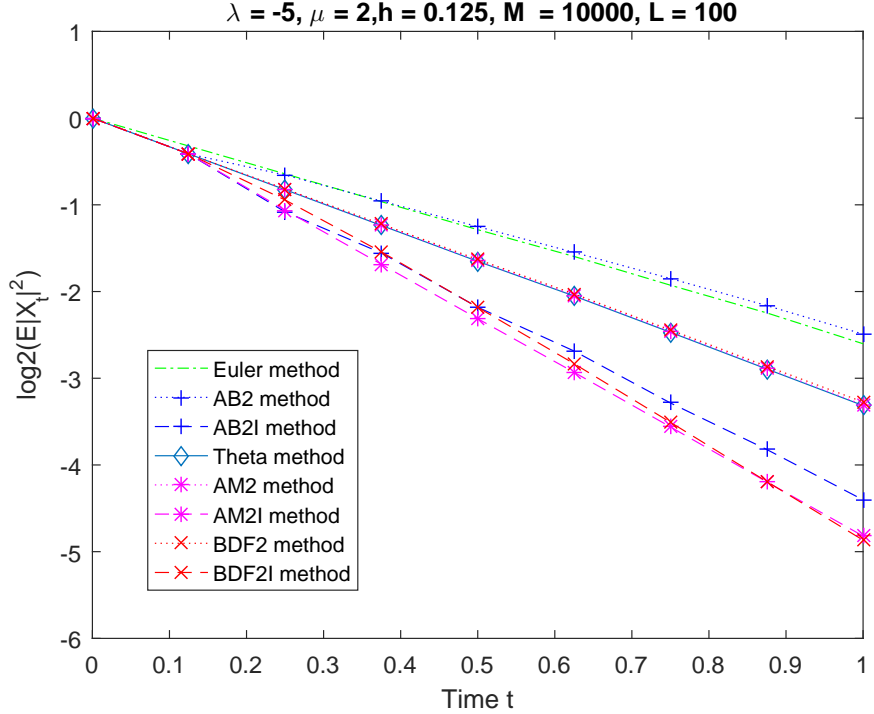}\label{STB-fig:exp1}
  \end{figure}

\begin{figure}[ht]
  \caption{Approximations of the $2^{nd}$ moment of the linear scalar test equation (\ref{STB-eq:sde.scalar}) in the interval $[0,20]$ for different two-step numerical  methods.}
  \centering
   \includegraphics[width=0.8\textwidth]{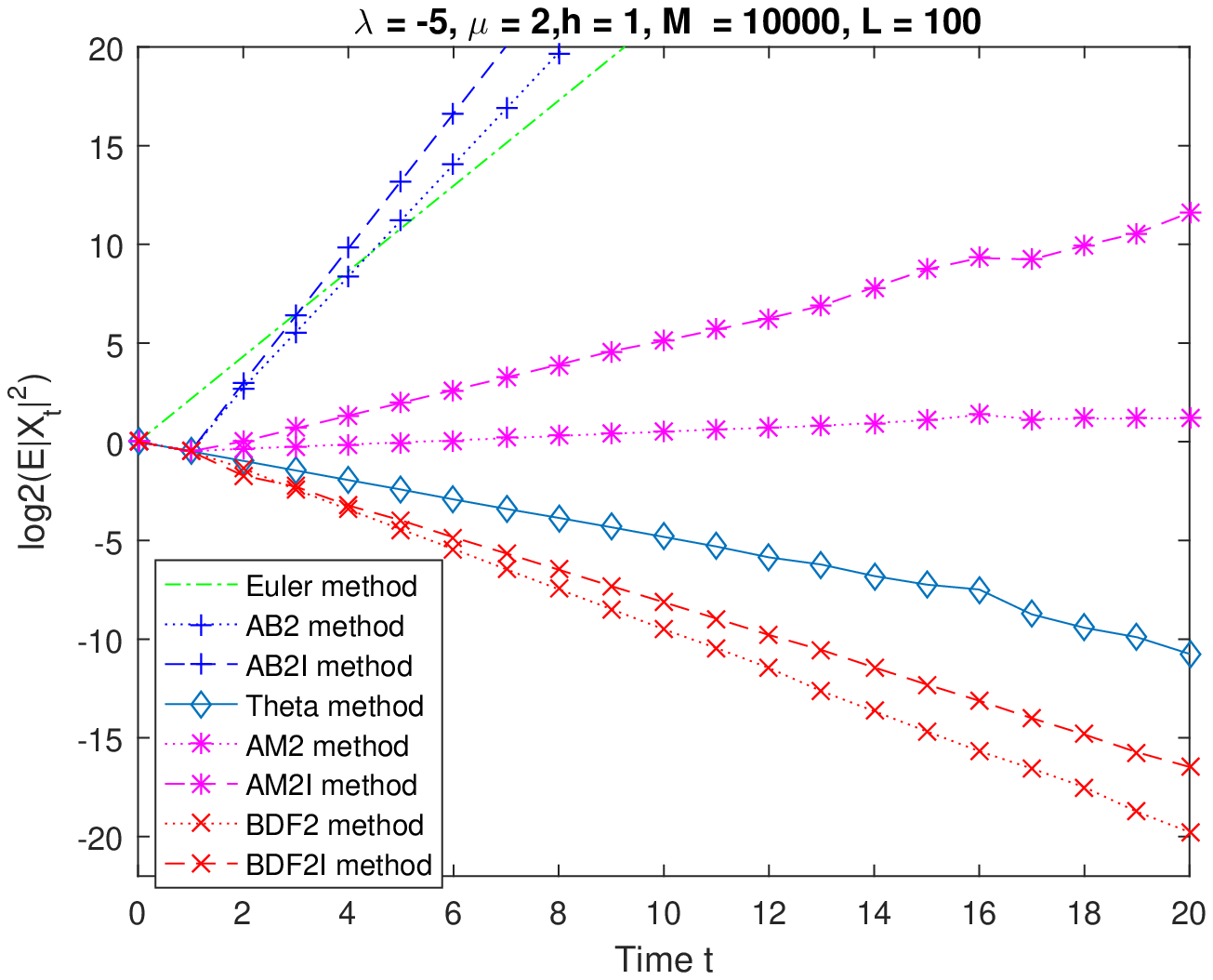}\label{STB-fig:exp20}
  \end{figure}

In this section we make some simple numerical experiments to complement the stability analysis presented in the previous section. We apply the AB2, AM2 and BDF2 two-step Maruyama schemes as well as their improved versions with constant step-size $h,$ to solve the equation
$$
dX_t = -5X_t + 2X_t dW_t, \qquad X_0=1.
$$ 
For the second initial condition in the two-step schemes we apply the $\theta-$Maruyama method which applied to the linear test equation (\ref{STB-eq:sde.scalar}) 
 reads
 $$
 X_{n+1}^{\theta EM} =\frac{1+ (1-\theta)\lam h + \mu\sqrt{h}\xi_n}{1 - \theta \lam h}X_n^{\theta EM} 
 $$
with $\theta=1/2$ and $n=0.$ We also implement the $\theta-$method (with $\theta=1/2$) and the Euler method (with $\theta=0$) for further comparison. We plot the obtained values  in a $log_2$-scale against time $t.$ The estimated mean-square norm of $X$ is point-wise estimated by each stochastic numerical method $X^{SM}$ in the following way,
$$
\sqrt{\bfE\left(X(t_i)^2\right)} \approx \left( \frac{1}{M L}\sum_{j=1}^{M}\sum_{k=1}^{L} \left(X^{SM}_{k,j}(t_i)\right)^2\right)^{1/2},
$$
where we have computed $M$ batches of $L$ simulation paths. The total number of paths in the experiments is $M\cdot L=10^6.$ For the first experiment, see Figure \ref{STB-fig:exp1},  we have applied all the methods with time-step size $h=1/8,$ so that they are all stable. The considered time interval is $[0, 1].$ In the second experiment,  see Figure \ref{STB-fig:exp20}, we integrate over $[0, 20]$  with $h=1.$ In this case the $AB2, AB2I$ and $AM2, AM2I$ methods are not stable as well as the forward Euler method. The $BDF2, BDF2I$ methods as well as the $\theta-$methods are asymptotically stable in the mean-square sense with the $BDF2$ performing better. 
Another remark we can make in the one-dimensional case is about the performance of the proposed improved methods with respect to their stability behavior, which seems to follow the rule that we do not gain more w.r.t to stability performance by using higher order schemes (multiple integrals for the approximation of the diffusion coefficient), as one can see from both Figures  \ref{STB-fig:exp1} and  \ref{STB-fig:exp20}. Nevertheless, the situation is different in more dimensions as shown n Section \ref{NSF:sec:system}.

\section{Linear system of SDEs and multi-dimensional noise.}\label{NSF:sec:system}
\setcounter{equation}{0}

Consider the $d$-system of linear test-equations with $m$-dimensional multiplicative noise
\beqq\label{STB-eq:sde.multisystem}
dX(t) = F X(t) dt +\sum_{r=1}^m G_r X(t) dW_r(t), \quad X(t_0)=X_0,
\eeqq
where $F,G$ are $d\times d$ real-valued matrices and assume w.l.o.g. that $X_0$ is non-random.

\subsection{Stability of two-step methods for linear system of SDEs driven by multi-dimensional noise.}

The two-step Maryuama method with an equidistant step-size $h$ and approximations $X_i=(X_{1,i},X_{2,i},\ldots,X_{n,i})^{T}$ of the solution of (\ref{STB-eq:sde.multisystem})  read
\beqq\label{STB-eq:twostepapp.system}
\sum_{j=0}^{2}\alpha_jX_{i-j} = h\sum_{j=0}^{2} \beta_j FX_{i-j} + \sum_{r=1}^m\sum_{j=1}^{2} \gamma_j G_rX_{i-j}\sqrt{h}\xi_{r, i-j}, \quad i=2, 3, \ldots,
\eeqq
and can be represented as 
\beqq\label{STB-eq:system.diff}
X_{i} = AX_{i-1} + CX_{i-2} + \sum_{r=1}^mB_rX_{i-1}\xi_{r, i-1} + \sum_{r=1}^mD_rX_{i-2}\xi_{r, i-2}, \, i=2, 3,\ldots,
\eeqq
where  the matrices  $A, C$ and $B_r, D_r$ are given by 
\beqq\label{STB-eq:twostepappcoef.systemAC}
A=(\alpha_0\bbi_d - h\beta_0F)^{-1}(-\alpha_1\bbi_d + h\beta_1F), \quad 
C=(\alpha_0\bbi_d - h\beta_0F)^{-1}(-\alpha_2\bbi_d + h\beta_2F)
\eeqq
\beqq\label{STB-eq:twostepappcoef.systemBD}
B_r=(\alpha_0\bbi_d - h\beta_0F)^{-1}\sqrt{h}\gamma_1G_r 
, \quad D_r=(\alpha_0\bbi_d - h\beta_0F)^{-1}\sqrt{h}\gamma_2G_r, \quad r=1, \ldots,m
\eeqq
and for the improved versions
\beqq\label{STB-eq:twostepappcoef.systemBDimpr}
B^*_r=B_r + (\alpha_0\bbi_d - h\beta_0F)^{-1} h^{3/2}(\gamma_1 + \eta_1)FG_r 
, \quad D^*_r=D_r + (\alpha_0\bbi_d - h\beta_0F)^{-1}h^{3/2}(\gamma_2 + \eta_2)FG_r,
\eeqq
for $r=1, \ldots,m.$

Here, the stability or transition  matrix $\bbs$ of the two-step method (\ref{STB-eq:system.diff})  applied to linear system of the form (\ref{STB-eq:sde.multisystem})  reads
\beqq\label{STB-eq:Smat.system}
\bbs= \begin{bmatrix}
    A\otimes A + \sum_{r=1}^mB_r\otimes B_r &  A\otimes C & C\otimes A & C\otimes C + \sum_{r=1}^mD_r\otimes D_r + R\\
    A\otimes \bbi_d & 0 & C\otimes \bbi_d &\sum_{r=1}^mD_r\otimes B_r \\
    \bbi_d\otimes A & \bbi_d\otimes C & 0 & \sum_{r=1}^mB_r\otimes D_r \\
    \bbi_{d^2} & 0 & 0  & 0
\end{bmatrix},
\eeqq
with $R=\sum_{r=1}^m(A\otimes D_r)(B_r\otimes \bbi_d) + \sum_{r=1}^m(D_r\otimes A)(\bbi_d\otimes B_r). $

A result of the type of  Theorem \ref{STB-thr:2step}, that is a conclusion about the 
asymptotically zero mean-square stability of the two-step method (\ref{STB-eq:twostepapp.system}) applied to the linear system (\ref{STB-eq:sde.multisystem}), is again related with equivalent conditions for the relation $\rho(\bbs)<1.$ Now the characteristic polynomial of the stability matrix $\bbs$ is of order $4d^2.$ The computational effort  of the Schur-Cohn test (SCJ) is  now bigger, but one can reduce it by halving the  dimensions of the matrix, whose positive-definite character needs to be checked at the expense of some easily checked   inequalities on linear combinations of the coefficients of the polynomial (c.f. \cite{anderson_jury:73}).

\subsection{A linear system of SDEs driven by a single noise term.}

\begin{figure}[ht]
  \caption{Approximations of the $MS$-norm of the linear system equation (\ref{STB-eq:sde.unisystem}) 
  in the interval $[0,3]$ for different two-step numerical  methods.}
  \centering
   \includegraphics[width=0.8\textwidth]{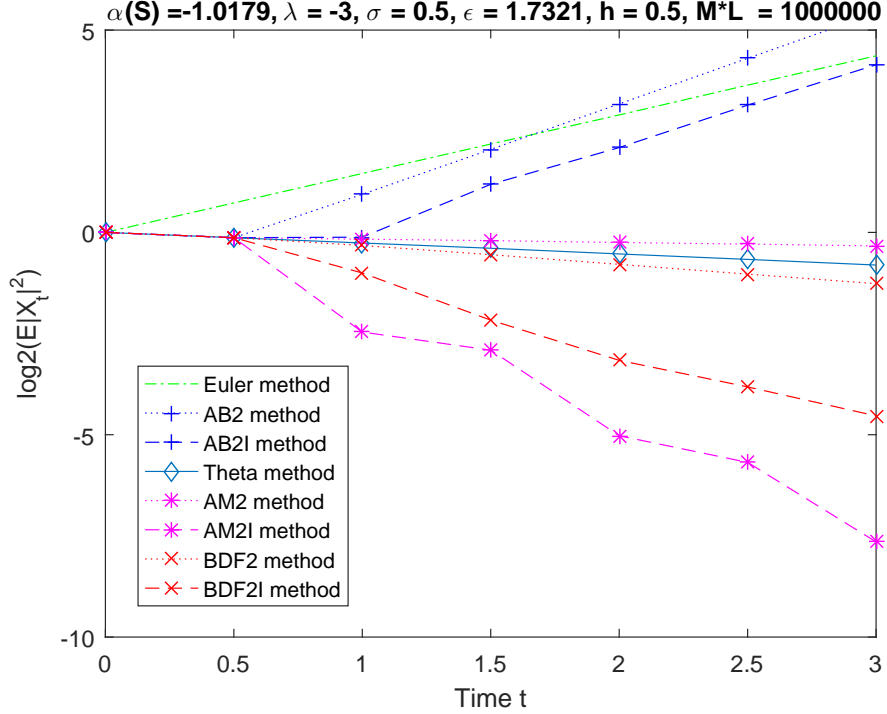}\label{STB-fig:exp3}
  \end{figure}

Consider the system of linear test-equations (\ref{STB-eq:sde.multisystem}) with $d=2, m=1$ and matrices $F, G$ of the following type
\beqq\label{STB-eq:systemF}
F= \begin{bmatrix}
    \lam & 0 \\
    0 &\lam
\end{bmatrix},   \quad G= \begin{bmatrix}
    \sigma & \ep \\
    \ep &\sigma
\end{bmatrix},
\eeqq
that is 
\beqq\label{STB-eq:sde.unisystem}
dX(t) = \begin{bmatrix}
    \lam & 0 \\
    0 &\lam
\end{bmatrix} X(t) dt +\begin{bmatrix}
    \sigma & \ep \\
    \ep &\sigma
\end{bmatrix}X(t) dW_1(t), \quad X(t_0)=X_0,
\eeqq
with a single noise term. The mean-square stability matrix for (\ref{STB-eq:sde.unisystem}) is
\beqq\label{STB-eq:systemS}
\bbs= \begin{bmatrix}
    2\lam + \sigma^2 & \sigma\ep & \ep\sigma & \ep^2\\
    \sigma\ep & 2\lam + \sigma^2 & \ep^2 &\ep\sigma \\
    \ep\sigma & \ep^2 & 2\lam + \sigma^2 & \sigma\ep \\
    \ep^2 & \ep\sigma & \sigma\ep  & 2\lam + \sigma^2
\end{bmatrix},
\eeqq
and the zero solution of (\ref{STB-eq:sde.unisystem}) is asymptotically MS-stable iff (cf. \cite[Lemma 4.1]{buckwar_sickenberger:12})
\beqq\label{STB-eq:systemuni.cond}
\lam + \frac{1}{2}\left( |\sigma| + |\ep|\right)^2<0.
\eeqq
Below we make a simple experiment implementing the two-step Maruyama methods 
\beqq\label{STB-eq:systemuni.diff}
X_{i} = AX_{i-1} + CX_{i-2} + BX_{i-1}\xi_{i-1} + DX_{i-2}\xi_{i-2}, \, i=2, 3,\ldots,
\eeqq
where in particular for the AB2/AB2I methods 
\beqq\label{STB-eq:twostepappcoef.systemACab2}
A=\bbi_2 + \frac{3}{2}hF, \quad C=-\frac{1}{2}hF,
\eeqq
\beqq\label{STB-eq:twostepappcoef.systemBDab2}
B=\sqrt{h}G, \quad B^*=\sqrt{h}G  + h^{3/2}FG, \quad D=0, \quad D^*=-\frac{1}{2}h^{3/2}FG, 
\eeqq
for the AM2/AM2I methods 
\beqq\label{STB-eq:twostepappcoef.systemACam2}
A=Q\left(\bbi_2 + \frac{8}{12}hF\right), \quad C=-\frac{1}{12}hQF,
\eeqq
\beqq\label{STB-eq:twostepappcoef.systemBDam2}
B=\sqrt{h}QG, \quad B^*=Q\left(\sqrt{h}G  + \frac{7}{12}h^{3/2}FG\right), \quad D=0, \quad D^*=-\frac{1}{12}h^{3/2}QFG, 
\eeqq
with $Q=(\bbi_2 - \frac{5}{12}hF)^{-1}$ and 
for the BDF2/BDF2I methods 
\beqq\label{STB-eq:twostepappcoef.systemACbdf2}
A=\frac{4}{3}Q, \quad C=-\frac{1}{3}Q,
\eeqq
\beqq\label{STB-eq:twostepappcoef.systemBDbdf2}
B=\sqrt{h}QG, \quad B^*=Q\left(\sqrt{h}G  + \frac{1}{3}h^{3/2}FG\right), \quad D=D^*=-\frac{1}{3}\sqrt{h}QG,  
\eeqq
with $Q=(\bbi_2 - \frac{2}{3}hF)^{-1}.$ We choose the values of $\lam, \sigma, \ep$ such that the spectral abscissa $\alpha(\bbs)$  of the mean-square stability matrix $\bbs$ is negative, that is $\alpha(\bbs)<0$,   and the spectral radius  $\rho(\bbs)<1$ or in other words such that the condition (\ref{STB-eq:systemuni.cond}) 
holds. In this case, see Figure \ref{STB-fig:exp3}, the improved versions AM2I and BDF2I are stable  whereas AM2 and BDF2 are not.

\subsection{A linear system of SDEs driven by two noise terms.}

Consider the system of linear test-equations (\ref{STB-eq:sde.multisystem}) with $d=2, m=2$ and matrices $F, G_1, G_2$ of the following type
\beqq\label{STB-eq:systemF2}
F= \begin{bmatrix}
    \lam & 0 \\
    0 &\lam
\end{bmatrix},   \quad G_1= \begin{bmatrix}
    \sigma & 0 \\
    0 &\sigma
\end{bmatrix}, \quad G_2= \begin{bmatrix}
    0 & -\ep \\
    \ep & 0
    \end{bmatrix},
\eeqq
that is 
\beqq\label{STB-eq:sde.unisystem2}
dX(t) = \begin{bmatrix}
    \lam & 0 \\
    0 &\lam
\end{bmatrix} X(t) dt +\begin{bmatrix}
    \sigma & 0 \\
    0 &\sigma
\end{bmatrix}X(t) dW_1(t) + \begin{bmatrix}
    0 & -\ep \\
    \ep & 0
\end{bmatrix}X(t) dW_2(t), \quad X(t_0)=X_0,
\eeqq
with two commutative noise terms. The mean-square stability matrix for (\ref{STB-eq:sde.unisystem2}) is
\beqq\label{STB-eq:systemS2}
\bbs= \begin{bmatrix}
    2\lam + \sigma^2 & 0 & 0 & \ep^2\\
    0 & 2\lam + \sigma^2 & -\ep^2 & 0\\
    0 & -\ep^2 & 2\lam + \sigma^2 & 0 \\
    \ep^2 & 0 & 0  & 2\lam + \sigma^2
\end{bmatrix},
\eeqq
and the zero solution of (\ref{STB-eq:sde.unisystem2}) is asymptotically MS-stable iff (cf. \cite[Lemma 4.1]{buckwar_sickenberger:12})
\beqq\label{STB-eq:systemuni2.cond}
\lam + \frac{1}{2}\left( \sigma^2 + \ep^2\right)<0.
\eeqq
Below we make a simple experiment implementing the two-step Maruyama methods 
\beqq\label{STB-eq:systemuni2.diff}
X_{i} = AX_{i-1} + CX_{i-2} + B_rX_{i-1}\xi_{r, i-1} + D_rX_{i-2}\xi_{r, i-2}, \, i=2, 3,\ldots,
\eeqq
where  for all methods $A$ and $B$ are as in (\ref{STB-eq:systemuni2.diff}) and  $B_r, D_r$ correspond now to the matrices $G_r;$ for instance for the AB2/AB2I methods we have
\beqq\label{STB-eq:twostepappcoef.system2BDab2}
B_r=\sqrt{h}G_r, \quad B_r^*=\sqrt{h}G_r  + h^{3/2}FG_r, \quad D_r=0, \quad D_r^*=-\frac{1}{2}h^{3/2}FG_r, \quad r=1,2.
\eeqq
 We choose the values of $\lam, \sigma, \ep$ in a way that the condition (\ref{STB-eq:systemuni2.cond}) holds and compute  the MS-norm of $X^{(1)},$ just as in  \cite{buckwar_sickenberger:12}, by
 $$
 \sqrt{\bfE(X_{t_i}^{(1)})^2} \approx \left(\frac{1}{ML}\sum_{j=1}^{M}\sum_{k=1}^{L}(X_{i,j,k}^{(1)}(\w))^2\right)^{1/2}.
 $$

\begin{figure}[ht]
  \caption{Approximations of the MS-norm of $X^{(1)}$ for the linear system equation (\ref{STB-eq:sde.unisystem}) in the interval $[0,3]$ with $h=1/2$ for different two-step numerical  methods.}
  \centering
   \includegraphics[width=0.8\textwidth]{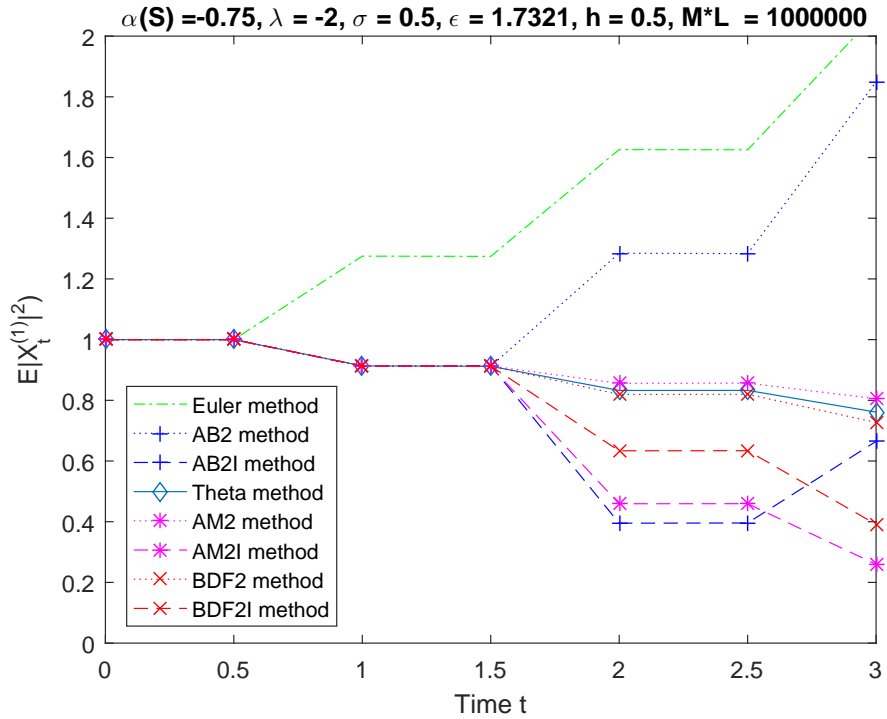}\label{STB-fig:exp4}
  \end{figure}

 In this case, see Figure \ref{STB-fig:exp4}, the improved versions AB2I, AM2I and BDF2I are stable  whereas AB2, AM2 and BDF2 are not.

\begin{figure}[ht]
  \caption{Approximations of the MS-norm of $X^{(1)}$ for the linear system equation (\ref{STB-eq:sde.unisystem}) in the interval $[0,3]$ with $h=1/4$ for different two-step numerical  methods.}
  \centering
   \includegraphics[width=0.8\textwidth]{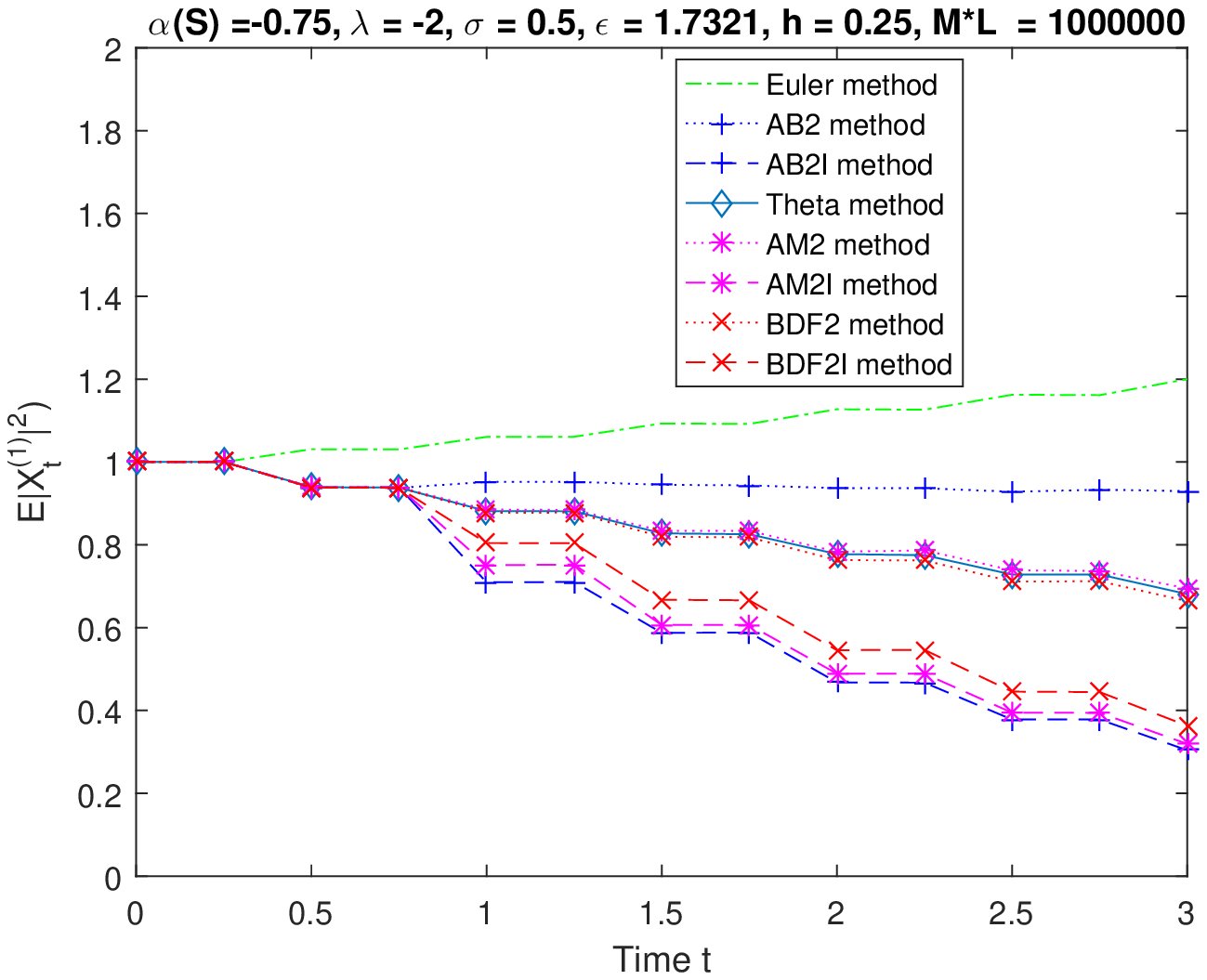}\label{STB-fig:exp5}
  \end{figure}

\begin{figure}[ht]
\centering
\begin{subfigure}{.45\textwidth}
   \includegraphics[width=1\textwidth]{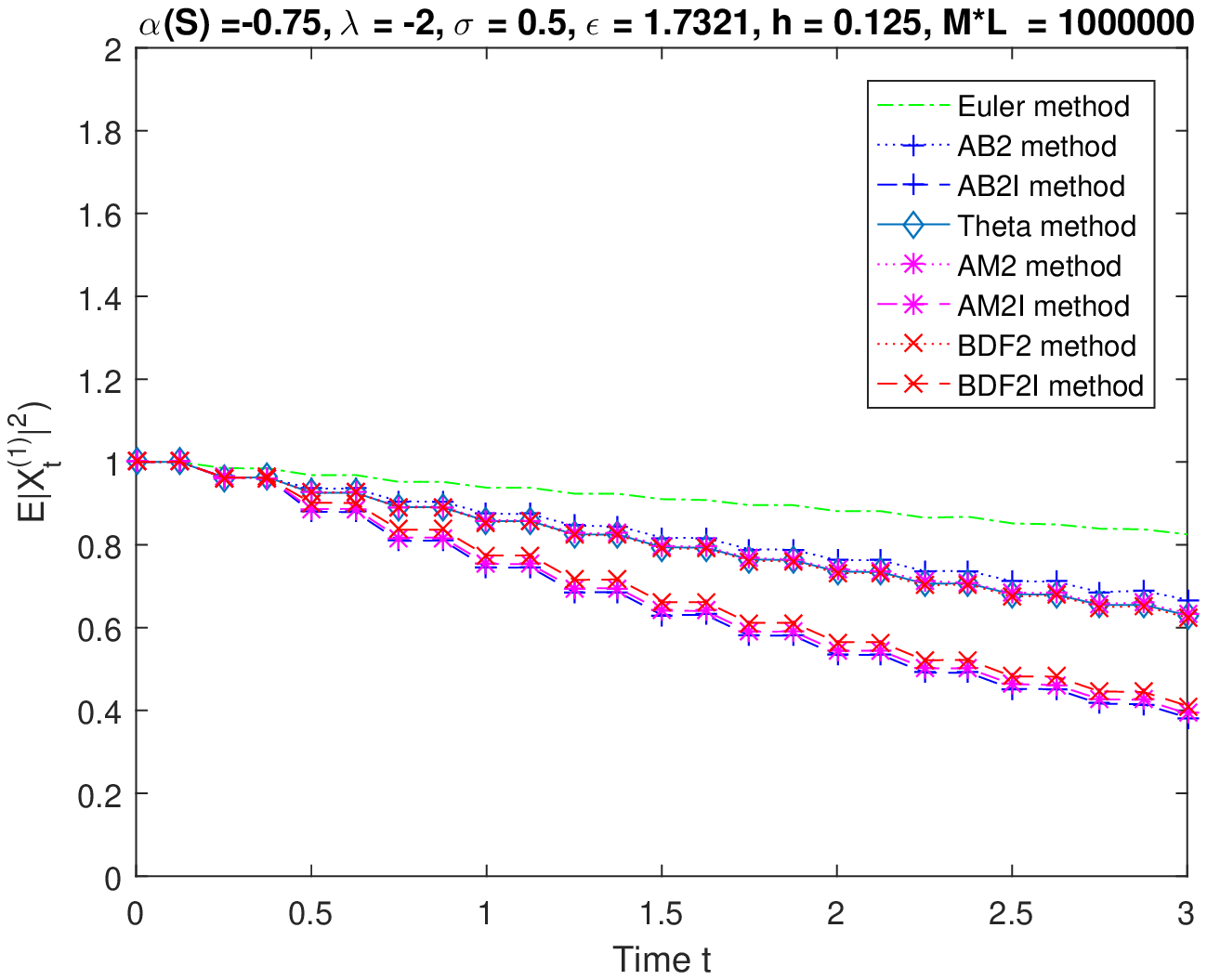}\label{STB-fig:exp6}
     \caption{Approximations of the MS-norm of $X^{(1)}$ for (\ref{STB-eq:sde.unisystem}) in the interval $[0,3]$ with $h=1/8.$}
  \end{subfigure}
  \begin{subfigure}{.45\textwidth}
  \centering
   \includegraphics[width=1\textwidth]{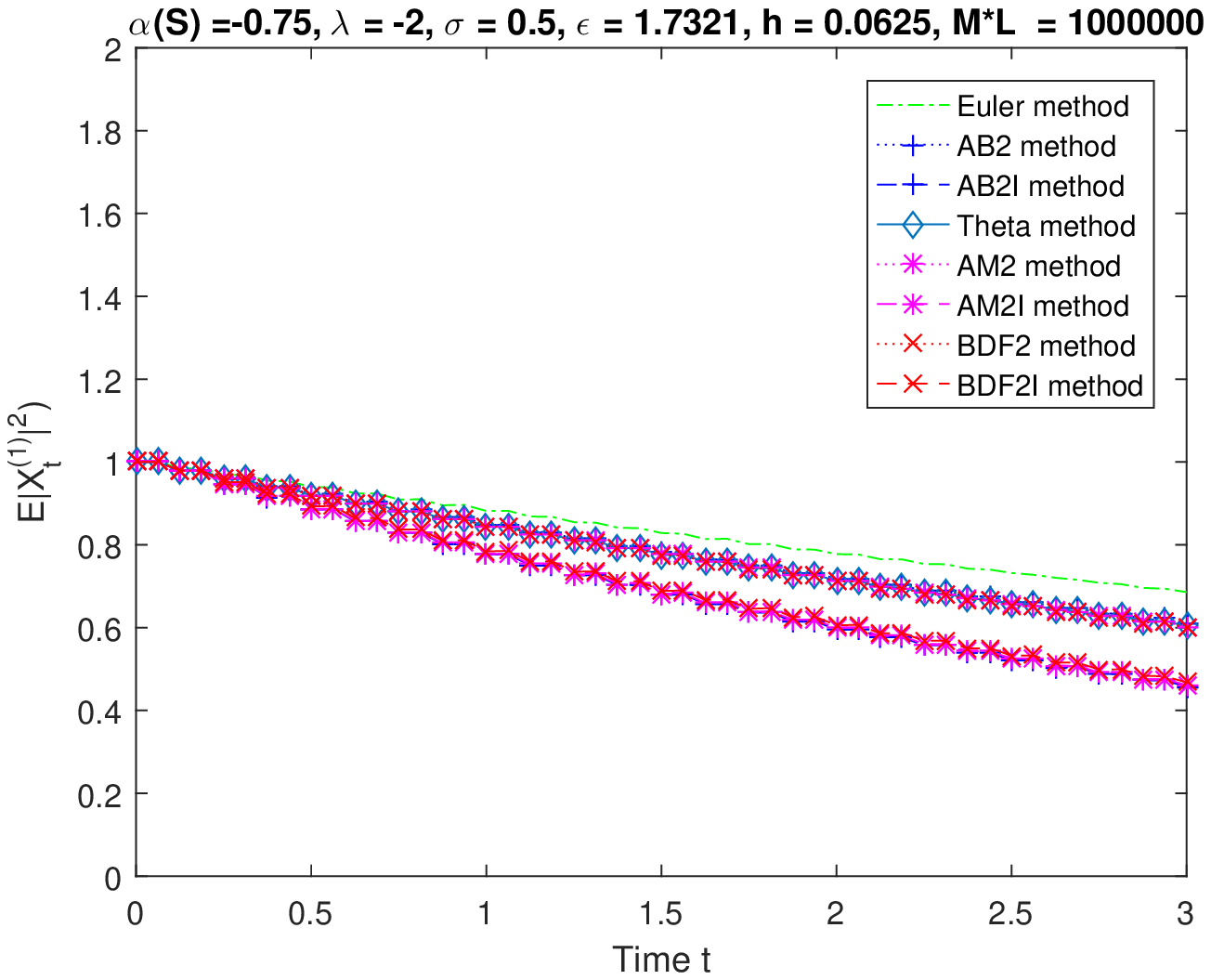}\label{STB-fig:exp7}
     \caption{Approximations of the MS-norm of $X^{(1)}$ for (\ref{STB-eq:sde.unisystem}) in the interval $[0,3]$ with $h=1/16.$}
\end{subfigure}
 \caption{Approximations of the MS-norm of $X^{(1)}$ for the linear system equation (\ref{STB-eq:sde.unisystem}) in the interval $[0,3]$ with $h=1/8, 1/16$ for different two-step numerical  methods.}\label{STB-fig:exp67}
  \end{figure}

\begin{figure}[ht]
\centering
\begin{subfigure}{.45\textwidth}
   \includegraphics[width=1\textwidth]{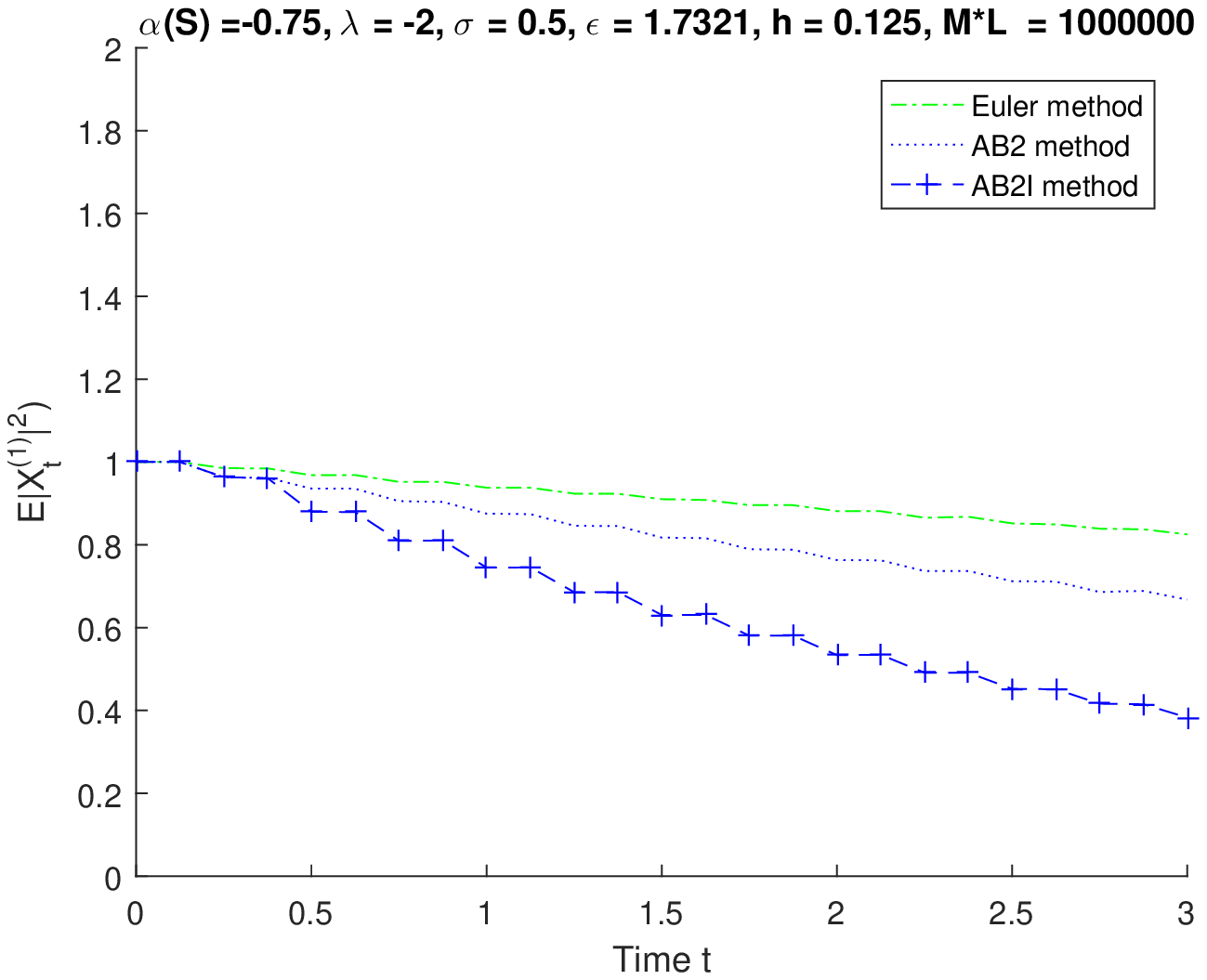}\label{STB-fig:exp6AB}
     \caption{Approximations of the MS-norm of $X^{(1)}$ for (\ref{STB-eq:sde.unisystem}) in the interval $[0,3]$ with $h=1/8$ for EM, AB2 and AB2I.}
  \end{subfigure}
  \begin{subfigure}{.45\textwidth}
  \centering
   \includegraphics[width=1\textwidth]{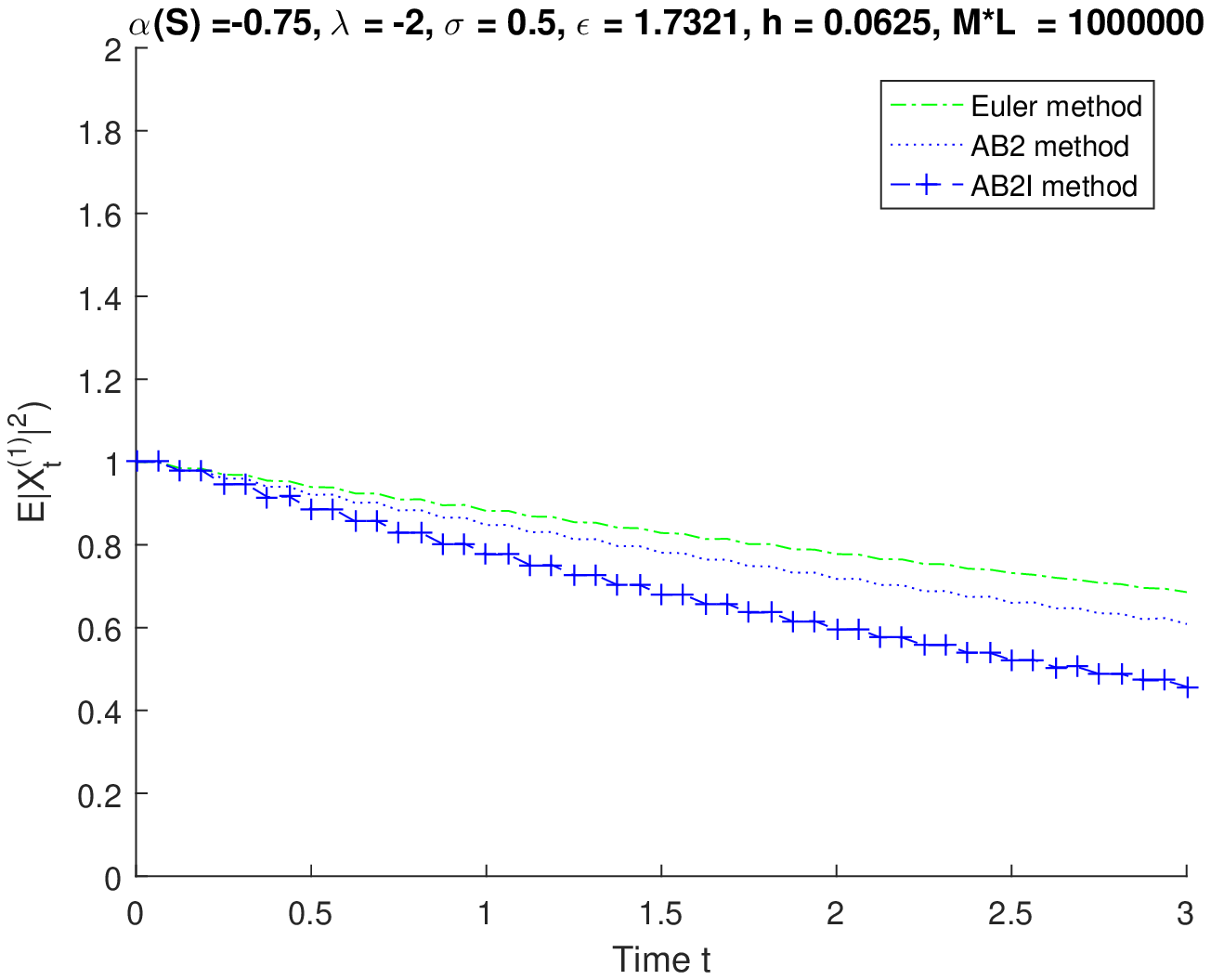}\label{STB-fig:exp7AB}
     \caption{Approximations of the MS-norm of $X^{(1)}$ for (\ref{STB-eq:sde.unisystem}) in the interval $[0,3]$ with $h=1/16$ for EM, AB2 and AB2I.}
\end{subfigure}
 \caption{Approximations of the MS-norm of $X^{(1)}$ for the linear system equation (\ref{STB-eq:sde.unisystem}) in the interval $[0,3]$ with $h=1/8, 1/16$ for EM, AB2 and AB2I.}\label{STB-fig:exp67AB}
  \end{figure}

\begin{figure}[ht]
\centering
\begin{subfigure}{.45\textwidth}
   \includegraphics[width=1\textwidth]{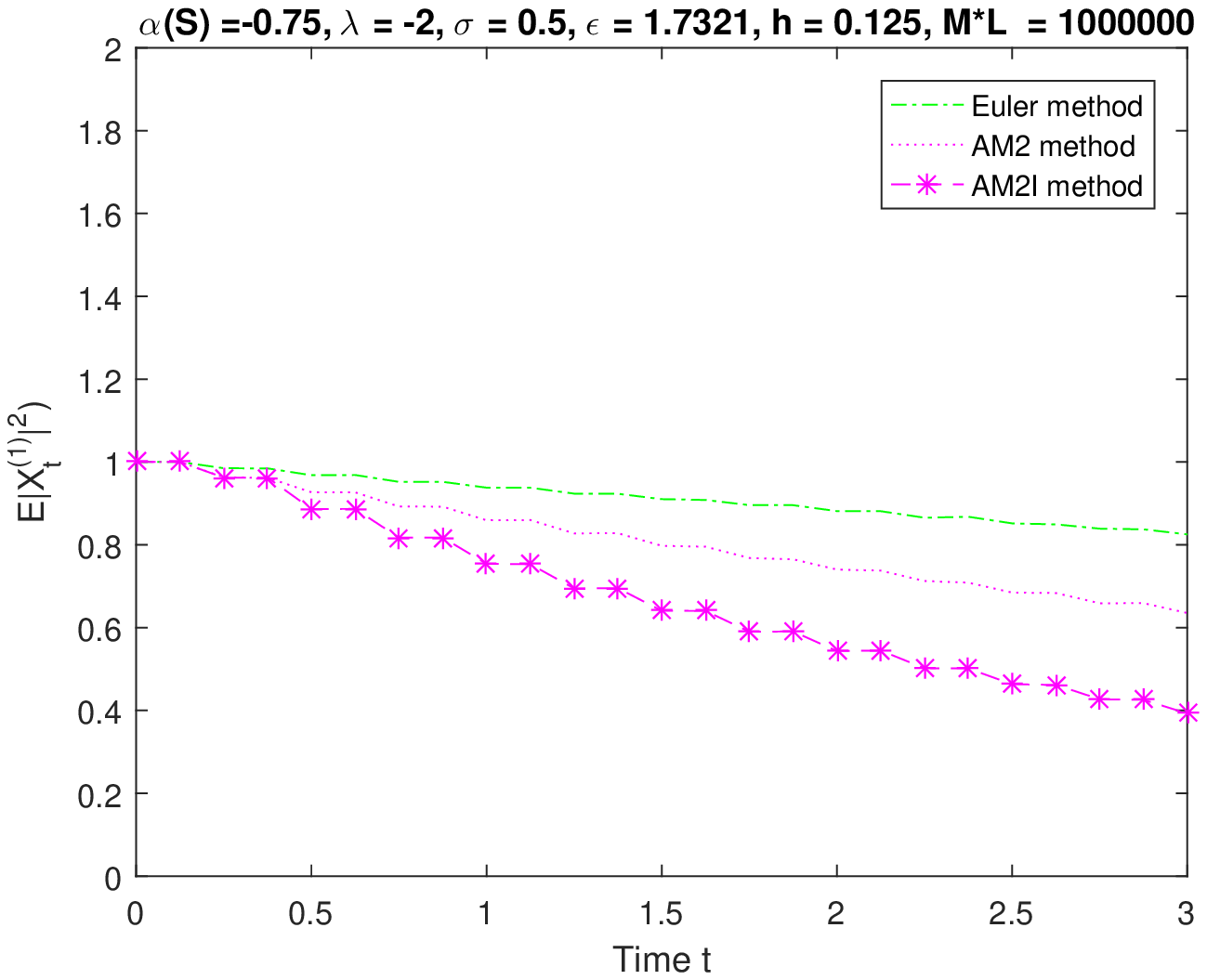}\label{STB-fig:exp6AM}
     \caption{Approximations of the MS-norm of $X^{(1)}$ for (\ref{STB-eq:sde.unisystem}) in the interval $[0,3]$ with $h=1/8$ for EM, AM2 and AM2I.}
  \end{subfigure}
  \begin{subfigure}{.45\textwidth}
  \centering
   \includegraphics[width=1\textwidth]{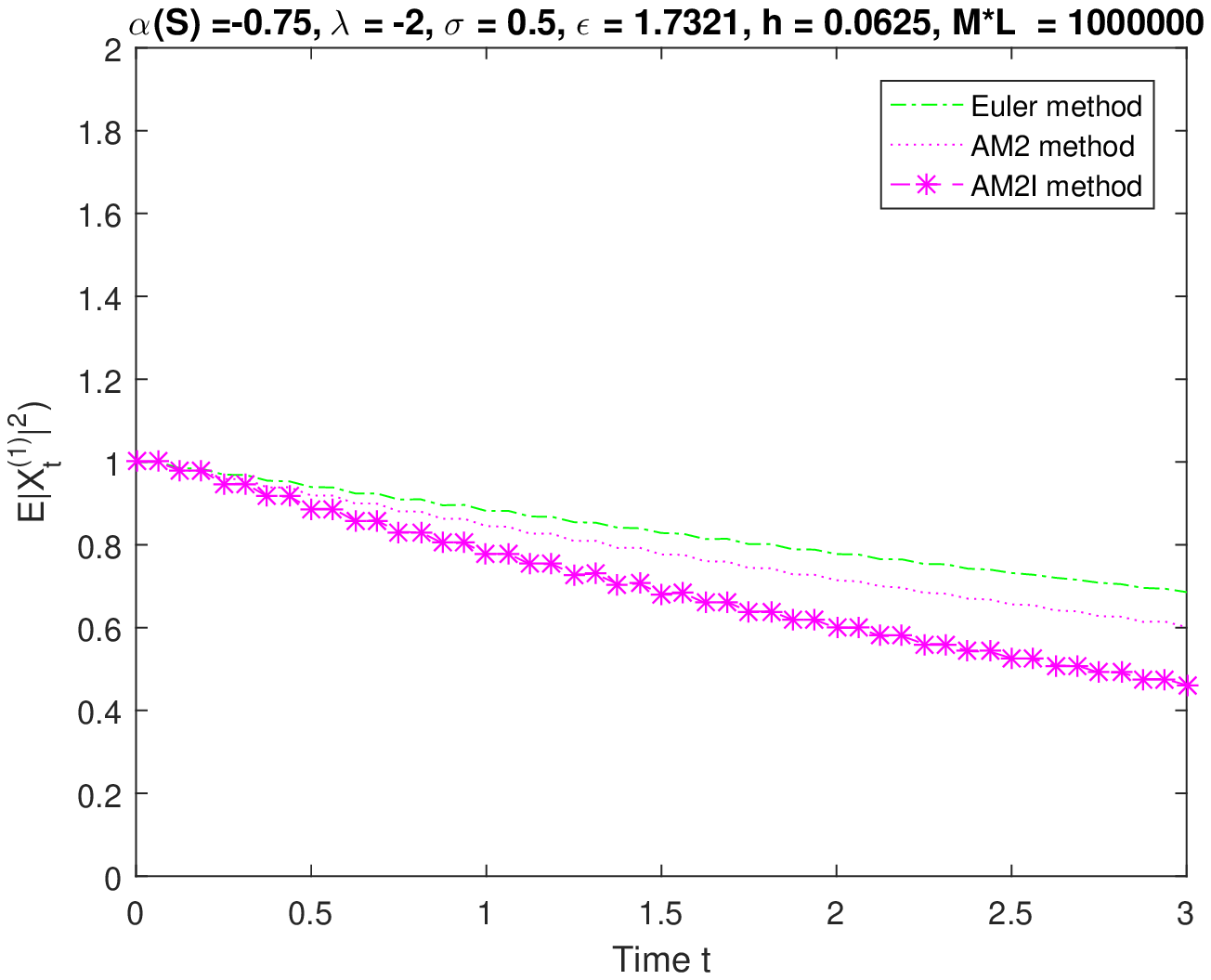}\label{STB-fig:exp7AM}
     \caption{Approximations of the MS-norm of $X^{(1)}$ for (\ref{STB-eq:sde.unisystem}) in the interval $[0,3]$ with $h=1/16$ for EM, AM2 and AM2I.}
\end{subfigure}
 \caption{Approximations of the MS-norm of $X^{(1)}$ for the linear system equation (\ref{STB-eq:sde.unisystem}) in the interval $[0,3]$ with $h=1/8, 1/16$ for EM, AM2 and AM2I.}\label{STB-fig:exp67AM}
  \end{figure}

\begin{figure}[ht]
\centering
\begin{subfigure}{.45\textwidth}
   \includegraphics[width=1\textwidth]{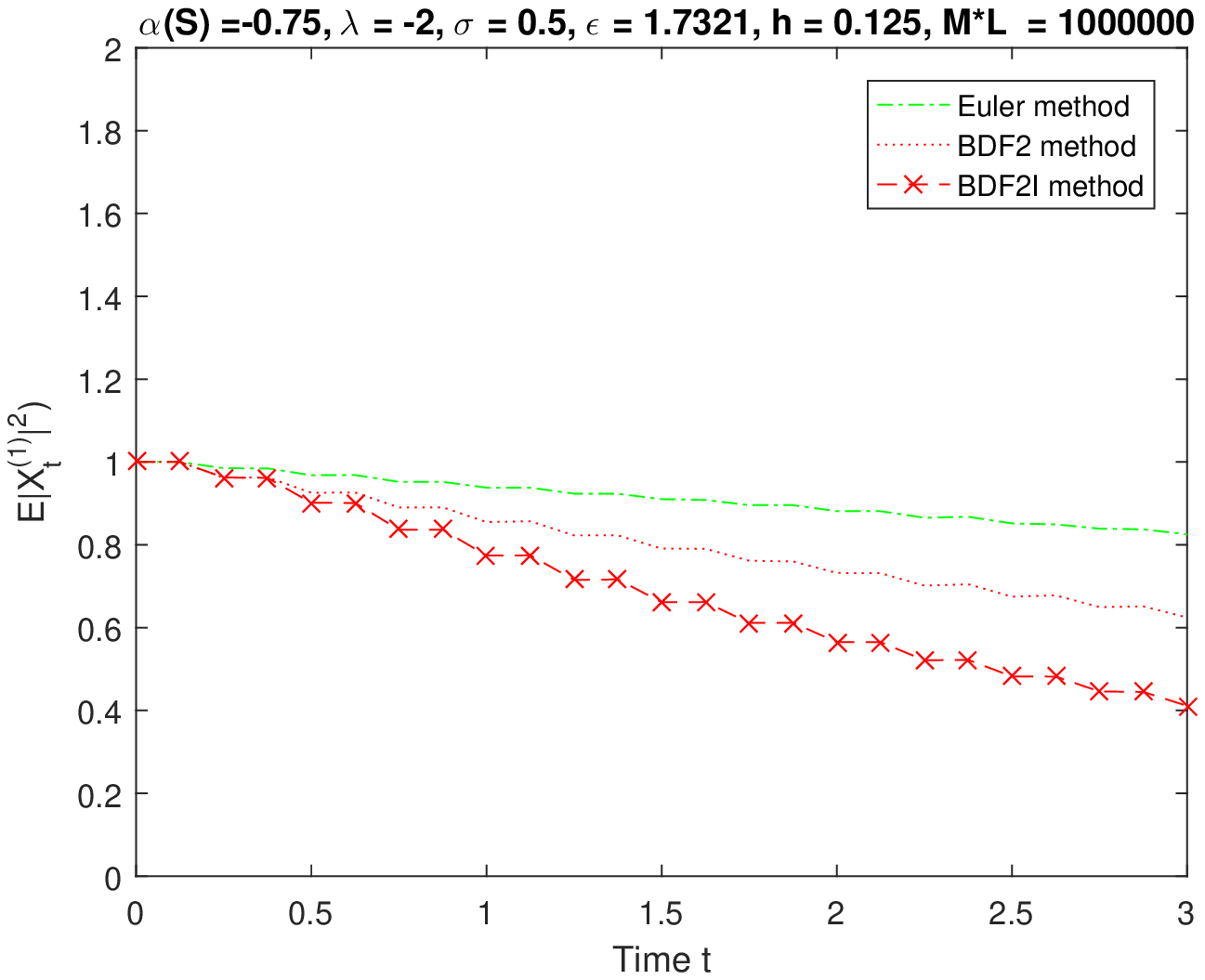}\label{STB-fig:exp6BDF}
     \caption{Approximations of the MS-norm of $X^{(1)}$ for (\ref{STB-eq:sde.unisystem}) in the interval $[0,3]$ with $h=1/8$ for EM, BDF2 and BDF2I.}
  \end{subfigure}
  \begin{subfigure}{.45\textwidth}
  \centering
   \includegraphics[width=1\textwidth]{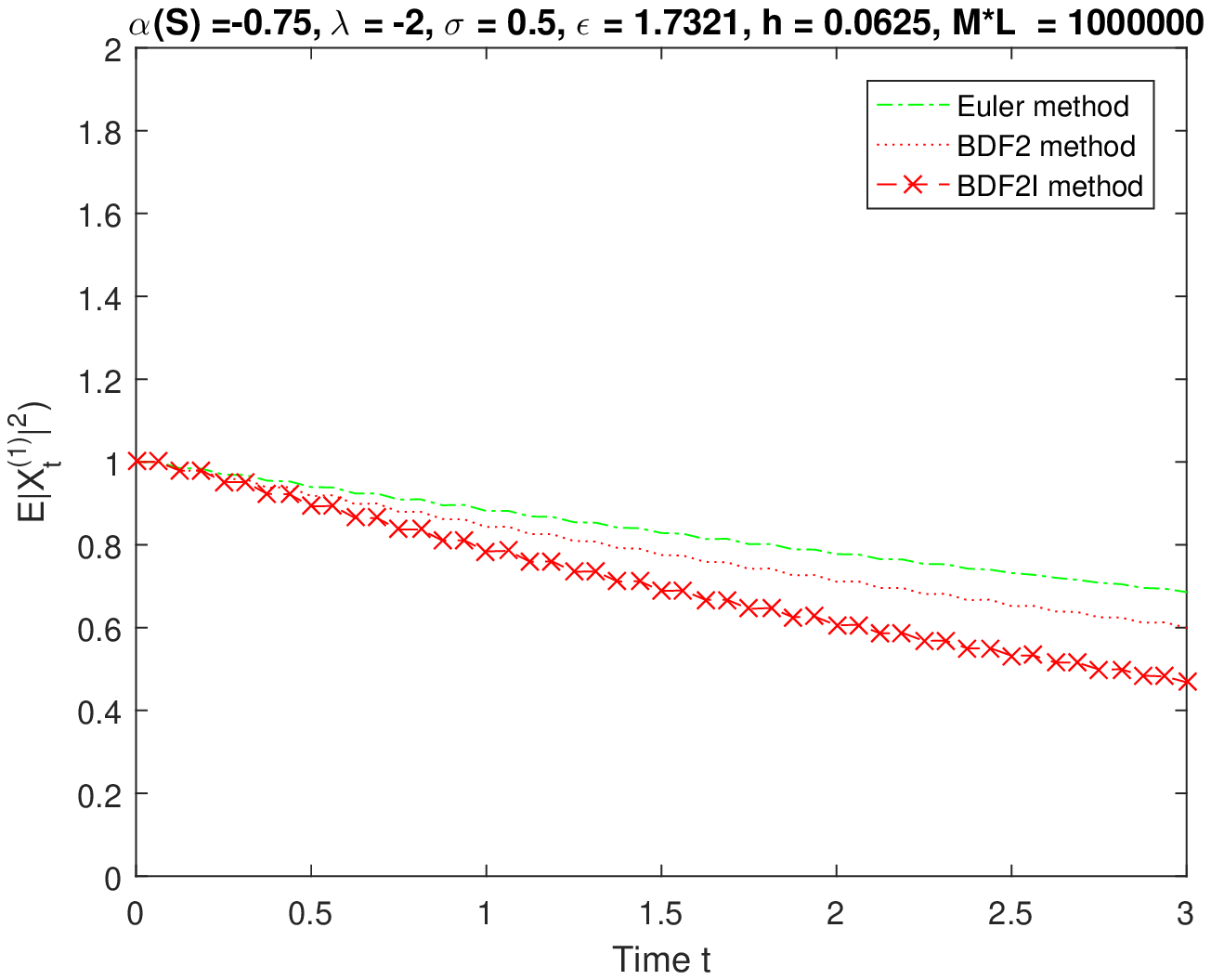}\label{STB-fig:exp7BDF}
     \caption{Approximations of the MS-norm of $X^{(1)}$ for (\ref{STB-eq:sde.unisystem}) in the interval $[0,3]$ with $h=1/16$ for EM, BDF2 and BDF2I.}
\end{subfigure}
 \caption{Approximations of the MS-norm of $X^{(1)}$ for the linear system equation (\ref{STB-eq:sde.unisystem}) in the interval $[0,3]$ with $h=1/8, 1/16$ for EM, BDF2 and BDF2I.}\label{STB-fig:exp67BDF}
  \end{figure}

Of course, by lowering the step-size the numerical methods become more stable. In the following, we sequentially halve the step-size and confirm the conjecture above. In all cases though, we conclude again that AB2, AM2 and BDF2 are less stable  than their improved counterparts, see Figures \ref{STB-fig:exp5},\ref{STB-fig:exp67} and for a clearer view Figures \ref{STB-fig:exp67AB},\ref{STB-fig:exp67AM} and \ref{STB-fig:exp67BDF}.

\subsection*{Acknowledgments}
The author would like to thank Prof. Evelyn Buckwar for fruitful discussions around the subject.
 It took place during the author's visit at the Institute for Stochastics, in Linz, Austria in the last four months of 2016.


\bibliographystyle{alpha}\baselineskip12pt 
\bibliography{stability}

\end{document}